\documentclass[acmsmall]{acmart}

\usepackage{amsmath}

\usepackage{xcolor}

\usepackage{algorithm}
\usepackage[noend]{algorithmic}

\usepackage{tikz}

\setcopyright{acmlicensed}
\copyrightyear{2026}
\acmYear{2026}
\acmDOI{XXXXXXX.XXXXXXX}

\acmJournal{TOMS}
\begin{document}

\title{Permutation-Avoiding FFT-Based Convolution}

\author{Nicolas Venkovic}
\email{venkovic@gmail.com, nicolas.venkovic@tum.de}
\author{Hartwig Anzt}
\email{hartwig.anzt@tum.de}
\affiliation{%
  \institution{Chair of Computational Mathematics, 
               School of CIT,
               TU Munich}
  \city{Heilbronn}
  \country{Germany}
}

\begin{abstract}
Fast Fourier Transform (FFT) libraries are widely used for evaluating discrete convolutions.
Most FFT implementations follow some variant of the Cooley-Tukey framework, in which the transform is decomposed into butterfly operations and index-reversal permutations.
While butterfly operations dominate the floating-point operation count, the memory access patterns induced by index-reversal permutations significantly degrade the FFT's arithmetic intensity.
When performing discrete convolution, the three sets of index-reversal permutations which occur in FFT-based implementations using Cooley-Tukey frameworks cancel out, thus paving the way to implementations free of any permutation.
To the best of our knowledge, such permutation-free variants of FFT-based discrete convolution are not commonly used in practice, making such kernels worth investigating.
Here, we look into such permutation-avoiding convolution procedures for multi-dimensional cases within a general radix Cooley-Tukey framework.
We perform numerical experiments to benchmark the algorithms presented against state-of-the-art FFT-based convolution implementations.
Our results suggest that developers of FFT libraries should consider supporting permutation-avoiding convolution kernels.
\end{abstract}





\maketitle

\section{Introduction}\label{sec:introduction}
Discrete convolution is a prevalent kernel of numerical simulations and data processing.
Example applications of discrete convolution are iterative solvers for classical physics simulation~\cite{moulinec1998numerical,zeman2010accelerating}, stationary Gaussian process simulation by circulant embedding~\cite{chan1997algorithm,dietrich1993fast,dietrich1997fast,wood1994simulation}, fast Poisson solvers~\cite{swarztrauber1977methods} as well as signal~\cite{oppenheim1999discrete} and image~\cite{gonzalez2009digital} processing.
One way to perform a discrete convolution is first to compute a forward discrete Fourier transform (DFT) of the data, do component-wise multiplication with the transformed filter, and then apply the inverse DFT.
In practice, this is achieved using highly optimized implementations of the forward and backward fast Fourier transform (FFT), e.g., FFTW~\cite{frigo1998fftw,frigo2005design},
Intel MKL~\cite{intelMKL}, etc.

The Cooley-Tukey framework, introduced in~\cite{cooley1965algorithm}, fundamentally transformed discrete Fourier transform computation by reducing the complexity from $\mathcal{O}(n^2)$ to $\mathcal{O}(n \log n)$ through a divide-and-conquer approach. 
While alternative FFT algorithms exist 
(prime factor algorithm~\cite{good1958interaction,thomas1963using}, 
Rader's algorithm~\cite{rader1968discrete}, 
Bluestein's algorithm~\cite{bluestein1970linear},
four-step FFT~\cite{gentleman1966fast,bailey1990ffts},
split-radix FFT~\cite{yavne1968economical,duhamel1984split}), the Cooley-Tukey framework remains the foundation of many production FFT libraries due to its broad applicability to composite sizes $n=r^t$, predictable performance, and implementation simplicity.
The Cooley-Tukey framework decomposes DFT computation into two distinct operations:
\begin{itemize}
\item \textbf{Butterfly operations}: Computational kernels with regular memory access patterns.
\item \textbf{Index-reversal permutations}: Data reorganization following index-reversal patterns, where index-reversal is a generalization of bit-reversal ordering to non-binary representations.
\end{itemize}
Here, our focus is on general radix Cooley-Tukey.
This choice is motivated by theoretical clarity in exposing computation-versus-data-movement trade-offs and practical relevance since insights directly translate to mixed-radix implementations used in production libraries.

While butterfly operations dominate floating-point operation count, index-reversal permutations create a fundamentally different performance bottleneck. 
Butterfly operations involve regular, predictable memory access patterns enabling vectorization and cache optimizations. 
In contrast, \emph{index-reversal permutations require scattered memory accesses following pseudo-random patterns, creating cache misses and preventing effective vectorization.}
Critically, \emph{index-reversal permutations have zero arithmetic intensity}, consisting purely of memory movement with no floating-point arithmetic, while butterfly operations achieve high arithmetic intensity through data reuse.
Index-reversal permutation complexity scales unfavorably with problem size. 
For 1D FFTs of size $n$, \emph{as $n$ increases beyond cache capacity, permutations operate on main memory with significantly higher latency}.

From the point of view of Cooley-Tukey, a standard FFT-based convolution requires three permutation applications per evaluation. 
How to best implement the associated data movements in the underlying FFT algorithms in relation to memory hierarchies has been the subject of many studies~(see, e.g., \cite{bailey1990ffts}) that incidentally led to optimal cache oblivious variants of FFT algorithms~\cite{frigo_cache-oblivious_2012}.
However, while alleviating the effects on runtime and energy consumption of unavoidable data movements in FFT algorithms is a well-documented practice, in the case of convolution, less attention has been spent on optimizing the data movement in convolutions.
Thus, as a means to achieve high performance implementations of convolution kernels, one may have to move away from mere composition of high-level kernels such as forward and backward transforms. 
Such an approach was for instance favored in~\cite{popovici_systematic_2021}.
In this work, we investigate permutation-avoiding FFT-based convolution of multi-dimensional data, that is, convolution free of any index-reversal permutations.
This approach, which was previously introduced as \emph{fast convolution} in \cite{van1992computational}, see also \cite{press_numerical_2007}, is not commonly used in practice.
However, in the current technological context, where the cost of moving data across memory keeps increasing relative to that of performing arithmetic operations, such kernels become increasingly relevant for achieving optimal performance.

This paper is organized as follows.
Section~\ref{sec:1d-problem} presents the 1D problem, Section~\ref{sec:2d-problem} extends to 2D, and Section~\ref{sec:multi-dimensional-problem} covers the general multi-dimensional case.
Section~\ref{sec:algorithms} presents algorithms and implementation details of the C code we developed in this study which is available at the following address:
\begin{center}
\url{https://github.com/venkovic/c-permutation-avoiding-convolution}.
\end{center}
Section~\ref{sec:numerical-experiments} provides experimental results, and Section~\ref{sec:discussion} discusses practical implications including applicability to other FFT frameworks.
Finally, the conclusions of this work are presented in Section~\ref{sec:conclusion}.

\section{1D problem}\label{sec:1d-problem}
The discrete Fourier transform (DFT) can be defined as
\begin{align}
\text{DFT}_n:
\mathbb{C}^n&\,\rightarrow \mathbb{C}^n\\
x&\,\mapsto F_nx\nonumber
\end{align}
where $F_n\in\mathbb{C}^{n\times n}$ is a Fourier matrix with components $(F_n)_{jk}=\omega^{j\cdot k}_n$ in which 
\begin{align}
\omega_{n}:=\exp(-2\pi i/n)
\end{align}
and where $i$ is the imaginary number.
The DFT being an invertible procedure, we have
\begin{align}
\text{DFT}_n^{-1}:
\mathbb{C}^n&\,\rightarrow \mathbb{C}^n\\
y&\,\mapsto F_n^{-1}y.\nonumber
\end{align}
Beyond the DFT, we are particularly interested in the function
\begin{align}\label{eq:convolution}
\mathcal{F}_h:\mathbb{C}^{n}&\rightarrow\mathbb{C}^n\\
x&\mapsto\text{DFT}_n^{-1}(\text{DFT}_n(h)\circ\text{DFT}_n(x))\nonumber
\end{align}
for some $h\in\mathbb{C}^n$, where $\circ$ denotes the Hadamard product.
We say that $\mathcal{F}_h(x)$ is the discrete convolution of $x$ with $h$, which is sometimes written as $x*h$.

\subsection{General radix Cooley-Tukey framework}\label{sec:cooley-tukey}
Several frameworks exist for formulating the FFT, the most common of which are based on variants of the approach introduced by James William Cooley and John Tukey~\cite{cooley1965algorithm}, such as general radix, mixed radix, and split-radix algorithms.
The Cooley-Tukey framework recursively decomposes a discrete Fourier transform (DFT) of composite size into smaller DFTs, achieving a computational complexity of $\mathcal{O}(n \log n)$.
Because it reduces the DFT into smaller subproblems, Cooley-Tukey can be combined with other DFT algorithms, such as Rader’s or Bluestein’s methods, to handle large prime sizes, or with the prime factor algorithm to exploit relatively prime factorizations more efficiently.
Moreover, the Cooley-Tukey method naturally admits a matrix-based formulation, which we present below, as it proves useful in describing our proposed method.
A good reference for what is presented below is the book of~\cite{van1992computational}.

If $n$ is a power of $r$, i.e., $n = r^t$ for some positive integer $t$, then we can resort to the general radix-$r$ Cooley-Tukey factorization~ of the Fourier matrix, i.e., we have
\begin{align}
F_{n}=&\,A_{r,n}P_{r,n}
\end{align}
in which $A_{r,n}\in\mathbb{C}^{n\times n}$ is a product of block butterfly matrices, i.e.,
\begin{align}
A_{r,n}=&\,B_{r,t,n}B_{r,t-1,n}\dots B_{r,2,n}B_{r,1,n}
\;\;\text{where}\;\;
t:=\log_r n
\end{align}
with block butterfly matrices given by
\begin{align}
B_{r,q,n}=I_{\frac{n}{r^q}}\otimes B_{r,r^q}=
\mathrm{diag}
(\underbrace{B_{r,r^q},\dots,B_{r,r^q}}_{\frac{n}{r^q}\;\mathrm{times}})
\;\;
\text{for }1\leq q\leq t,
\end{align}
where $\otimes$ denotes the Kronecker product, and $P_{r,n}\in\mathbb{R}^{n\times n}$ is a symmetric index-reversal permutation matrix, i.e., $P_{r,n}=P_{r,n}^T=P_{r,n}^{-1}$.

For any $n$ which is a power of $r$, the $n\times n$ butterfly matrix denoted by $B_{r,n}$ is given by
\begin{align}\label{eq:25}
B_{r,n}=
(F_r\otimes I_{n/r})\cdot
\text{diag}
\left(I_{n/r},\Omega_{r,n/r},\Omega_{r,n/r}^2,\dots,\Omega_{r,n/r}^{r-1}\right)
\end{align}
where 
\begin{align}
\Omega_{r,n}=\text{diag}
\left(1,\omega_{rn},\omega_{rn}^2,\dots,\omega_{rn}^{n-1}\right).
\end{align}
For example, in the radix-2 case, it is well-known that Eq.~\eqref{eq:25} simplifies to 
\begin{align}
B_{2,n}=
\begin{bmatrix}
I_{n/2}&\Omega_{n/2}\\
I_{n/2}&-\Omega_{n/2}
\end{bmatrix}.
\end{align}

The index-reversal permutation matrix denoted by $P_{r,n}\in\mathbb{R}^{n\times n}$ is given by
\begin{align}
P_{r,n}=R_{r,t,n}R_{r,t-1,n}\dots R_{r,2,n}R_{r,1,n}
\end{align}
where each permutation matrix $R_{r,q,n}$ is of the form
\begin{align}
R_{r,q,n}=
I_{r^{q-1}}\otimes R_{r,n/r^{q-1}}=
\mathrm{diag}(\underbrace{R_{r,n/r^{q-1}},\dots,R_{r,n/r^{q-1}}}_{r^{q-1}\;\mathrm{times}})
\end{align}
in which $R_{r,n/{r^{q-1}}}$ is a modulo $r$ sort permutation matrix.
That is, for any $n$ which is a power of $r$, $R_{r,n}$ is such that
\begin{align}\label{eq:r-ary-sorting}
R_{r,n}x=
[x_1,x_{r+1},x_{2r+1},\dots,x_{n-r+1},x_2,x_{r+2},x_{2r+2},\dots,x_{n-r+2},\dots,x_r,x_{2r},x_{3r},\dots,x_{n}]^T
\end{align}
for all $x\in\mathbb{C}^{n}$. 
This operation groups elements by their residue modulo-$r$, effectively performing a radix-$r$ digit-based reordering.
The modulo-$r$ sorting generalizes the even-odd sorting from the radix-2 case in the context of zero-based indexing.

Note that, by symmetry of the Fourier matrix, we have
\begin{align}
F_{n}=F_{n}^T=(A_{r,n}P_{r,n})^T=P^T_{r,n}A_{r,n}^T=P_{r,n}A_{r,n}^T.
\end{align}
It is also well known (e.g., see~\cite{van1992computational}) that
\begin{align}
F_{n}^{-1}=
\overline{F_{n}}/n=
\overline{A_{r,n}P_{r,n}}/n=
\overline{A_{r,n}}P_{r,n}/n.
\end{align}

\subsection{Permutation-avoiding 1D convolution}
Using the factorizations $F_n=P_{r,n}A_{r,n}^T$ and $F_n^{-1}=\overline{A_{r,n}}P_{r,n}/n$ in the matrix form of Eq.~\eqref{eq:convolution}, we obtain
\begin{align}
\mathcal{F}_h(x)=
F_{n}^{-1}((F_nh)\circ (F_{n}x))=
\overline{A_{r,n}}P_{r,n}((P_{r,n}A_{r,n}^Th)\circ (P_{r,n}A_{r,n}^Tx))/n
\end{align}
where the permutation $P_{r,n}$ in $(P_{r,n}x)\circ(P_{r,n}y)$ can be factored out as $P_{r,n}(x\circ y)$, so that
\begin{align}\label{eq:34}
\mathcal{F}_h(x)=
\overline{A_{r,n}}P_{r,n}P_{r,n}((A_{r,n}^Th)\circ (A_{r,n}^Tx))/n=
\overline{A_{r,n}}((A_{r,n}^Th)\circ (A_{r,n}^Tx))/n.
\end{align}
in which $P_{r,n}P_{r,n}=P_{r,n}P_{r,n}^{-1}=I_n$.
Therefore, we have:
\begin{align}\label{eq:35}
\mathcal{F}_h(x)=
\overline{A_{r,n}}((A_{r,n}^Th)\circ (A_{r,n}^Tx))/n.
\end{align}
Eq.~\eqref{eq:35} is what we refer to as the permutation-avoiding 1D convolution of $x$ with $h$, i.e., $x*h$.
In~\cite{van1992computational}, this is referred to as \textit{fast convolution}.
The evaluation of $\mathcal{F}_h$ after Eq.~\eqref{eq:35} is index-reversal-free, whereas a standard FFT-based evaluation of the convolution requires \emph{three applications of the index-reversal permutation per evaluation of} $\mathcal{F}_h$.

\section{2D problem}\label{sec:2d-problem}
The 2D DFT can be defined as
\begin{align}
\text{DFT}_{n_1\times n_2}:
\mathbb{C}^{n_1\times n_2}&\,\rightarrow \mathbb{C}^{n_1\times n_2}\\
X&\,\mapsto F_{n_1}XF_{n_2}\nonumber
\end{align}
with an inverse given by
\begin{align}
\text{DFT}_{n_1\times n_2}^{-1}:
\mathbb{C}^{n_1\times n_2}&\,\rightarrow \mathbb{C}^{n_1\times n_2}\\
Y&\,\mapsto F_{n_1}^{-1}YF_{n_2}^{-1}.\nonumber
\end{align}
Analogously to the 1D problem, we are interested in the function
\begin{align}\label{eq:convolution-2d}
\mathcal{F}_H:\mathbb{C}^{n_1\times n_2}&\rightarrow\mathbb{C}^{n_1\times n_2}\\
X&\mapsto\text{DFT}_{n_1\times n_2}^{-1}(\text{DFT}_{n_1\times n_2}(H)\circ\text{DFT}_{n_1\times n_2}(X))\nonumber
\end{align}
for some $H\in\mathbb{C}^{n_1\times n_2}$.
$\mathcal{F}_H(X)$ is the discrete convolution of $X$ with $H$ denoted by $X*H$.
In matrix form, we have
\begin{align}\label{eq:convolution-2d-mtx}
\mathcal{F}_H(X)=&\,F_{n_1}^{-1}((F_{n_1}HF_{n_2})\circ (F_{n_1}XF_{n_2}))F_{n_2}^{-1}.
\end{align}
Let us now introduce row-major column vectorization, namely
\begin{align}
\text{vec}:\mathbb{C}^{n_1\times n_2}&\,\rightarrow\mathbb{C}^{n_1n_2}\\
X&\,\mapsto\begin{bmatrix}(X_{1,1:n_2})^T\\(X_{2,1:n_2})^T\\\vdots\\(X_{n_1,1:n_2})^T\end{bmatrix}.\nonumber
\end{align}
Then, for all $A\in\mathbb{C}^{n_1\times n_1}$, $X\in\mathbb{C}^{n_1\times n_2}$ and $B\in\mathbb{C}^{n_2\times n_2}$, we have 
\begin{align}
\text{vec}(AXB)=(A\otimes B^T)\cdot\text{vec}(X)
\end{align}
where, due to the symmetry of Fourier matrices, we have
\begin{align}\label{eq:45}
\text{vec}(F_{n_1}XF_{n_2})
=(F_{n_1}\otimes F_{n_2}^T)\cdot\text{vec}(X)
=(F_{n_1}\otimes F_{n_2})\cdot\text{vec}(X)
\end{align}
and, similarly, for the inverse transforms:
\begin{align}\label{eq:50}
\text{vec}(F_{n_1}^{-1}YF_{n_2}^{-1})
=(F_{n_1}^{-1}\otimes F_{n_2}^{-T})\cdot\text{vec}(Y)
=(F_{n_1}^{-1}\otimes F_{n_2}^{-1})\cdot\text{vec}(Y).
\end{align}
Using Eq.~\eqref{eq:50} along with Eq.~\eqref{eq:convolution-2d-mtx}, we get
\begin{align}
\text{vec}(\mathcal{F}_H(X))=
(F_{n_1}^{-1}\otimes F_{n_2}^{-1})\cdot\text{vec}\left((F_{n_1}H F_{n_2})\circ(F_{n_1}X F_{n_2})\right).
\end{align}
Note that, for all $X$ and $Y$, we have $\text{vec}(X\circ Y)=\text{vec}(X)\circ\text{vec}(Y)$.
Therefore, we have:
\begin{align}
\text{vec}(\mathcal{F}_H(X))=
(F_{n_1}^{-1}\otimes F_{n_2}^{-1})\cdot\left(\text{vec}(F_{n_1}H F_{n_2})\circ\text{vec}(F_{n_1}X F_{n_2})\right).
\end{align}
in which using Eq.~\eqref{eq:45} leads to
\begin{align}\label{eq:55}
\text{vec}(\mathcal{F}_H(X))=
(F_{n_1}^{-1}\otimes F_{n_2}^{-1})
\cdot
\left(((F_{n_1}\otimes F_{n_2})\cdot\text{vec}(H))
\circ
((F_{n_1}\otimes F_{n_2})\cdot\text{vec}(X))\right).
\end{align}

\subsection{Permutation-avoiding 2D convolution}
Using $F_n=P_{r,n}A_{r,n}^T$ along with the fact that $(AB)\otimes(CD)=(A\otimes C)\cdot(B\otimes D)$, we get
\begin{align}\label{eq:60}
F_{n_1}\otimes F_{n_2}=
(P_{r,n_1}A^T_{r,n_1})\otimes(P_{r,n_2}A^T_{r,n_2})=
(P_{r,n_1}\otimes P_{r,n_2})\cdot(A^T_{r,n_1}\otimes A^T_{r,n_2}).
\end{align} 
Similarly, using $F_n^{-1}=\overline{A_{r,n}}P_{r,n}/n$, we get
\begin{align}\label{eq:65}
F_{n_1}^{-1}\otimes F_{n_2}^{-1}=
(\overline{A_{r,n_1}}P_{r,n_1})\otimes(\overline{A_{r,n_2}}P_{r,n_2})/(n_1n_2)=
(\overline{A_{r,n_1}}\otimes \overline{A_{r,n_2}})\cdot(P_{r,n_1}\otimes P_{r,n_2})/(n_1n_2)
\end{align}
so that using Eqs.~\eqref{eq:60} into Eq.~\eqref{eq:55} leads to
\begin{align}\label{eq:70}
\text{vec}&\,(\mathcal{F}_H(X))=\\
&\hspace{-.45cm}(F_{n_1}^{-1}\otimes F_{n_2}^{-1})
\cdot
\left[
\left((P_{r,n_1}\otimes P_{r,n_2})\cdot(A^T_{r,n_1}\otimes A^T_{r,n_2})\cdot\text{vec}(H)\right)
\circ
\left((P_{r,n_1}\otimes P_{r,n_2})\cdot(A^T_{r,n_1}\otimes A^T_{r,n_2})\cdot\text{vec}(X)\right)
\right]
\nonumber
\end{align}
in which the permutation matrix $(P_{r,n_1}\otimes P_{r,n_2})$ may be factored out as follows:
\begin{align}
\text{vec}(\mathcal{F}_H(X))=
(F_{n_1}^{-1}\otimes F_{n_2}^{-1})
\cdot
(P_{r,n_1}\otimes P_{r,n_2})\cdot
\left[
\left((A^T_{r,n_1}\otimes A^T_{r,n_2})\cdot\text{vec}(H)\right)
\circ
\left((A^T_{r,n_1}\otimes A^T_{r,n_2})\cdot\text{vec}(X)\right)
\right].
\end{align}
Using Eq.~\eqref{eq:65} along with the fact that $(P_{r,n_1}\otimes P_{r,n_2})\cdot(P_{r,n_1}\otimes P_{r,n_2})=I_{n_1n_2}$, we obtain:
\begin{align}
\text{vec}(\mathcal{F}_H(X))=
(\overline{A_{r,n_1}}\otimes \overline{A_{r,n_2}})
\cdot
\left[
\left((A^T_{r,n_1}\otimes A^T_{r,n_2})\cdot\text{vec}(H)\right)
\circ
\left.
\left((A^T_{r,n_1}\otimes A^T_{r,n_2})\cdot\text{vec}(X)\right)
\right]
\right/(n_1n_2),
\end{align}
so that:
\begin{align}\label{eq:80}
\mathcal{F}_H(X)=
\text{vec}^{-1}
\left.
\left(
(\overline{A_{r,n_1}}\otimes \overline{A_{r,n_2}})
\cdot
\left[
\left((A^T_{r,n_1}\otimes A^T_{r,n_2})\cdot\text{vec}(H)\right)
\circ
\left((A^T_{r,n_1}\otimes A^T_{r,n_2})\cdot\text{vec}(X)\right)
\right]\right)
\right/
(n_1n_2).
\end{align}
Similarly as for the 1D case, Eq.~\eqref{eq:80} allows for the evaluation of the convolution of $X$ with $H$, i.e., $X*H$, without the need to apply index-reversal permutations.
On the other hand, a naive application of Eqs.~\eqref{eq:convolution-2d}-\eqref{eq:convolution-2d-mtx} requires \emph{three applications of the permutation matrix} $P_{r,n_1}\otimes P_{r,n_2}$ \emph{for each evaluation of} $\mathcal{F}_H$.

\section{Multi-dimensional problem}\label{sec:multi-dimensional-problem}
Consider the $d$-dimensional DFT applied to a $d$-way array $\mathscr{X}\in\mathbb{C}^{n_1 \times \cdots \times n_d}$. 
The forward DFT is defined as
\begin{align}
\text{DFT}_{n_1\times \dots\times n_d}:
\mathbb{C}^{n_1\times\dots\times n_d}&\,\rightarrow \mathbb{C}^{n_1\times \dots\times n_d}\\
\mathscr{X}&\,\mapsto \text{vec}^{-1}((F_{n_1}\otimes\dots\otimes F_{n_d})\cdot\text{vec}(\mathscr{X}))\nonumber
\end{align}
and its inverse is given by
\begin{align}
\text{DFT}^{-1}_{n_1\times \dots\times n_d}:
\mathbb{C}^{n_1\times\dots\times n_d}&\,\rightarrow \mathbb{C}^{n_1\times \dots\times n_d}\\
\mathscr{Y}&\,\mapsto \text{vec}^{-1}((F^{-1}_{n_1}\otimes\dots\otimes F^{-1}_{n_d})\cdot\text{vec}(\mathscr{Y}))\nonumber
\end{align}
in which the invertible row-major column vectorization operator for a $d$-way array $\mathscr{X} \in \mathbb{C}^{n_1\times\cdots\times n_d}$ is given by
\begin{align}
\text{vec}:\mathbb{C}^{n_1 \times \dots \times n_d}&\,\rightarrow \mathbb{C}^{n_1 \cdots n_d}\\
\mathscr{X} &\,\mapsto 
[\mathscr{X}_{1\dots 11,1:n_d}^T,
\dots,
\mathscr{X}_{1\dots 1n_{d-1},1:n_d}^T,
\mathscr{X}_{1\dots 21,1:n_d}^T,
\dots,
\mathscr{X}_{1\dots 2n_{d-1},1:n_d}^T,
\dots,
\mathscr{X}_{n_1\dots n_{d-1}n_d}^T]^T. \nonumber
\end{align}
That is, the entries of $\mathscr{X}$ are stacked by increasing index in the last mode, then the before-last, and so on, with the last index varying the fastest.
Note that the $\otimes$ operator in those definitions of the multi-dimensional DFT and its inverse is still a Kronecker product, not a tensor product, i.e., $F_{n_1}\otimes\dots\otimes F_{n_d}$ is a $n\times n$ matrix with $n=n_1\cdots n_d$, that is, a 2-way array, irrespective of $d$.

To perform a convolution of $\mathscr{X}\in\mathbb{C}^{n_1\times\dots\times n_d}$ with $\mathscr{H}\in\mathbb{C}^{n_1\times\dots\times n_d}$, we define
\begin{align}\label{eq:convolution-multi-dimensional}
\mathcal{F}_\mathscr{H}:\mathbb{C}^{n_1\times\dots\times n_d}&\rightarrow\mathbb{C}^{n_1\times\dots\times n_d}\\
\mathscr{X}&\mapsto\text{DFT}_{n_1\times\dots\times n_d}^{-1}
\left(\text{DFT}_{n_1\times\dots\times n_d}(\mathscr{H})
\circ
\text{DFT}_{n_1\times\dots\times n_d}(\mathscr{X})\right).\nonumber
\end{align}
Upon vectorizing $\mathscr{X}$, which stacks the $d$-way array into a column vector, and using Kronecker product $\otimes$ identities, we obtain:
\begin{align}\label{eq:90}
\text{vec}(\mathcal{F}_\mathscr{H}(\mathscr{X})) =
(F_{n_1}^{-1}\otimes\dots\otimes F_{n_d}^{-1})
\cdot
\left[
\left((F_{n_1}\otimes\dots\otimes F_{n_d})\cdot\text{vec}(\mathscr{H})\right)
\circ
\left((F_{n_1}\otimes\dots\otimes F_{n_d})\cdot\text{vec}(\mathscr{X})\right)
\right].
\end{align}

\subsection{Permutation-avoiding multi-dimensional convolution}
Using the factorizations of Fourier matrix as $F_{n} = P_{r,n} A_{r,n}^T$ and $F_n^{-1}=\overline{A_{r,n}}P_{r,n}/n$, and factoring out permutation matrices as we did in previous sections, we are then able to recast Eq.~\eqref{eq:90} as follows:
\begin{align}
\text{vec}(\mathcal{F}_\mathscr{H}(\mathscr{X}))=\\
&\hspace{-1.9cm}(\overline{A_{r,n_1}}\otimes\dots\otimes\overline{A_{r,n_d}})
\cdot
\left.
\left[
((A_{r,n_1}^T\otimes\dots\otimes A_{r,n_d}^T)\cdot\text{vec}(\mathscr{H}))
\circ
((A_{r,n_1}^T\otimes\dots\otimes A_{r,n_d}^T)\cdot\text{vec}(\mathscr{X}))
\right]
\right/(n_1\cdots n_d).\nonumber
\end{align}
Once again, this formulation allows for the evaluation of the convolution $\mathcal{F}_{\mathscr{H}}(\mathscr{X})$ without the need to carry out index-reversal permutations.

\section{Algorithms}\label{sec:algorithms}
Although our presentation of the Cooley-Tukey framework is matrix-based, its traditional implementation, and what we implemented in \url{https://github.com/venkovic/c-permutation-avoiding-convolution}, are not.
Similarly, although the basic idea of Cooley-Tukey algorithms is recursive, most implementations, including ours, rearrange the algorithm to avoid explicit recursion and optimize memory access patterns.
A more detailed description of the fundamental kernels used throughout this study is provided in Section~\ref{app}.

As is representative of common practice---our implementation included---the performance of FFT algorithms is evaluated assuming precomputation of twiddle factors.
In general, a greater radix leads to smaller numbers of floating-point operations (FLOP).
In particular, for Cooley-Tukey algorithms, we have:
\begin{itemize}
\item Radix-2: $5\log_2 n$ FLOP.
\item Radix-4: $4.25\log_2 n$ FLOP.
\item Radix-8: $4.08\log_2 n$ FLOP.
\end{itemize}

However, the number of FLOP only draws an incomplete picture of actual performance. 
In practice, memory access complexity also carries a significant effect on real-world performance.
For a complete FFT evaluation, the memory access pattern consists of: 
(i) butterfly operations, which exhibit regular stride patterns with high spatial and temporal locality; 
(ii) index-reversal permutations, which create pseudo-random access patterns with poor cache behavior; and 
(iii) twiddle factor access, which is generally regular and cacheable when precomputed.

To better understand algorithm performance on modern processors, the roofline model provides a framework that relates arithmetic intensity to achievable performance.
Arithmetic intensity, defined as the ratio of floating-point operations to memory traffic, determines whether an algorithm is compute-bound or memory-bound on a given architecture.
While butterfly operations achieve high arithmetic intensity through the reuse of loaded data elements, permutations have zero arithmetic intensity—every memory access serves only to move data without computation.

The key insight is that index-reversal permutations, despite not involving any floating-point arithmetic, can consume a significant fraction of total execution time due to their memory-bound nature.
For larger problem sizes that exceed cache capacity, our FFT kernels operate in memory-bound regimes where permutation costs become the dominant performance bottleneck.
This memory-bound behavior explains why theoretical FLOP reductions from higher radices do not necessarily translate to performance improvements when permutation overhead overwhelms the computational savings.

Our implementation of the Cooley-Tukey FFT is meticulously decomposed into modular and transparent components: a precomputed index-reversal permutation, a precomputed twiddle factor array, and efficient in-place, split-complex, unordered FFT (or butterfly) kernels. 
This deliberate separation allows us to not only evaluate the overall FFT cost but also to isolate and quantify the contribution of the index-reversal permutation step and the main butterfly computation independently.
This modular approach provides several advantages.
First, we can precisely measure the cost of permutations versus butterfly operations.
Second, changes to permutation handling can be evaluated in isolation.
Third, memory access costs can be studied independently of computational costs.
Lastly, the impact of different cache hierarchies and memory systems can be assessed.
Crucially, this structured approach mirrors the practical use-case modeled by FFTW when called with the \texttt{FFTW\_ESTIMATE} flag: that is, a scenario where no advanced architecture-specific tuning and optimizations are applied beyond basic algorithmic design.

All code is compiled using gcc with the -march=native compiler option on. 
While this allows for the use of all instruction sets supported by the machine on which the code is compiled, our butterfly implementations do not make explicit use of any vectorization instructions, even though our use of the -O3 option may trigger some of those at compilation. 
A more finely tuned implementation of the butterfly routines that would make use of specific instruction sets could increase the relative runtime attributed to index permutations, making permutation-avoiding kernels more effective.

The algorithms presented below represent different trade-offs between computational efficiency and memory access patterns, providing a foundation for understanding when permutation-avoiding strategies offer the greatest benefit.
\begin{itemize}
\item Algo.~\ref{alg:index-reversal-permutation} (\textbf{Index-reversal permutation}) isolates the pure memory movement cost, providing a baseline for understanding the overhead that our permutation-avoiding approach eliminates.
\item  Algos.~\ref{alg:uFFT1D}-\ref{alg:uFFT3D} (\textbf{Unordered forward FFTs}) perform only butterfly operations without any permutations, achieving optimal arithmetic intensity but producing results in index-reversed order.
Note that, due to the symmetry of Fourier matrices, the Cooley-Tukey framework can be implemented in either of two ways: first, the data is permuted before a standard butterfly operation $\text{uFFT}_{r,n}(x) =A_{r,n}x$ is applied; second, a transposed set of butterfly operations is applied through $\text{uFFT2}_{r,n}(x) =A^T_{r,n}x$, prior to permuting the data.
We implement both types of unordered forward transforms.
\item  Algo.~\ref{alg:forward-fft} (\textbf{Forward FFT}) represents the standard approach, combining both permutation and butterfly costs.
\item  Algo.~\ref{alg:uiFFT1D}-\ref{alg:uiFFT3D} (\textbf{Unordered backward FFTs}) perform only conjugate butterfly operations without any permutations, achieving optimal arithmetic intensity but producing results in index-reversed order.
\item  Algo.~\ref{alg:backward-fft} (\textbf{Backward FFT}) represents the standard approach, combining both permutation and conjugate butterfly costs.
\item  Algo.~\ref{alg:fft-based-convolution} (\textbf{Standard FFT-based convolution}) demonstrates the triple permutation cost in traditional convolution implementations.
\item  Algo.~\ref{alg:permutation-avoiding-convolution} (\textbf{Permutation-avoiding convolution}) shows the approach which eliminates permutations while maintaining correctness.
\end{itemize}
For small problems that fit entirely in cache, permutations may represent a small fraction of the total FFT runtime. 
However, as problem size grows beyond cache capacity, permutation costs should represent a greater part of the execution time due to their memory-bound nature, particularly for higher-dimensional problems where the memory access patterns become increasingly irregular.
This scaling behavior creates a compelling case for permutation-avoiding approaches, especially in scenarios involving large problem sizes where cache effects become pronounced, and memory-constrained systems where bandwidth is the primary bottleneck.

\begin{algorithm}[ht]
\caption{Index-reversal permutation, $\text{permute}_{r,n_1\times\cdots\times n_d}(x)$}
\label{alg:index-reversal-permutation}
\begin{algorithmic}[1]
\REQUIRE{$x:=\text{vec}(\mathscr{X})\in\mathbb{C}^{n_1\cdots n_d},\mathscr{X}\in\mathbb{C}^{n_1\times\dots\times n_d},r\text{ s.t. }\exists\; t_1,\dots,t_d\text{ s.t. }n_1=r^{t_1},\dots,n_d=r^{t_d}$}
\RETURN{$(P_{r,n_1}\otimes\dots\otimes P_{r,n_d})\,x$}
\end{algorithmic}
\end{algorithm}

\begin{algorithm}[ht]
\caption{Unordered 1D forward FFT, $\text{uFFT(2)}_{r,n}(x)$}
\label{alg:uFFT1D}
\begin{algorithmic}[1]
\REQUIRE{$x\in\mathbb{C}^{n}, r\text{ s.t. }\exists\; t\text{ s.t. }n=r^{t}$}
\IF{$\text{uFFT}_{r,n}$}
\RETURN{$A_{r,n}x$}
\COMMENT{See Algos.~\ref{alg:butterfly-kernel-radix-2}, \ref{alg:butterfly-kernel-radix-4} and \ref{alg:butterfly-kernel-radix-8}}
\ELSIF{$\text{uFFT2}_{r,n}$}
\RETURN{$A_{r,n}^Tx$}
\COMMENT{See Algos.~\ref{alg:transposed-butterfly-kernel-radix-2}, \ref{alg:transposed-butterfly-kernel-radix-4} and \ref{alg:transposed-butterfly-kernel-radix-8}}
\ENDIF
\end{algorithmic}
\end{algorithm}

\begin{algorithm}[ht]
\caption{Unordered 2D forward FFT, $\text{uFFT(2)}_{r,n_1\times n_2}(x)$}
\label{alg:uFFT2D}
\begin{algorithmic}[1]
\REQUIRE{$x:=\text{vec}(X)\in\mathbb{C}^{n_1n_2},X\in\mathbb{C}^{n_1\times n_2},\;r\text{ s.t. }\exists\; t_1,t_2\;\text{ s.t. }\;n_1=r^{t_1},n_2=r^{t_2}$}
\ENSURE{$(A_{r,n_1}\otimes A_{r,n_2})\,x$ or $(A^T_{r,n_1}\otimes A^T_{r,n_2})\,x$}
\FOR{$i_1=1,\dots,n_1$}
\STATE{$i_{2_0}:=(i_1-1)n_2$}
\STATE{$[x_{i_{2_0}+1},\dots,x_{i_{2_0}+n_2}]:=
\text{uFFT(2)}_{r,n_2}([x_{i_{2_0}+1},\dots,x_{i_{2_0}+n_2}]^T)^T$}
\ENDFOR
\FOR{$i_2=1,\dots,n_2$}
\STATE{$i_{1_0}:=-1$}
\STATE{$z:=[x_{(i_{1_0}+1)n_2+i_2},\dots,x_{(i_{1_0}+n_1)n_2+i_2}]$}
\STATE{$z:=\text{uFFT(2)}_{r,n_1}(z^T)^T$}
\STATE{$[x_{(i_{1_0}+1)n_2+i_2},\dots,x_{(i_{1_0}+n_1)n_2+i_2}]:=z$}
\ENDFOR
\RETURN{$x$}
\end{algorithmic}
\end{algorithm}

\begin{algorithm}[ht]
\caption{Unordered 3D forward FFT, $\text{uFFT(2)}_{r,n_1\times n_2\times n_3}(x)$}
\label{alg:uFFT3D}
\begin{algorithmic}[1]
\REQUIRE{$x:=\text{vec}(\mathscr{X})\in\mathbb{C}^{n_1n_2n_3},\;\mathscr{X}\in\mathbb{C}^{n_1\times n_2\times n_3},\;r\text{ s.t. }\exists\; t_1,t_2,t_3\;\text{ s.t. }\;n_1=r^{t_1},n_2=r^{t_2},n_3=r^{t_3}$}
\ENSURE{$(A_{r,n_1}\otimes A_{r,n_2}\otimes A_{r,n_3})\,x$ or $(A^T_{r,n_1}\otimes A^T_{r,n_2}\otimes A_{r,n_3}^T)\,x$}
\FOR{$i_1=1,\dots,n_1$}
\FOR{$i_2=1,\dots,n_2$}
\STATE{$i_{3_0}:=((i_1-1)n_2+i_2-1)n_3$}
\STATE{$[x_{i_{3_0}+1},\dots,x_{i_{3_0}+n_3}]:=
\text{uFFT(2)}_{r,n_3}([x_{i_{3_0}+1},\dots,x_{i_{3_0}+n_3}]^T)^T$}
\ENDFOR
\ENDFOR
\FOR{$i_1=1,\dots,n_1$}
\FOR{$i_3=1,\dots,n_3$}
\STATE{$i_{2_0}:=(i_1-1)n_2-1$}
\STATE{$z:=[x_{(i_{2_0}+1)n_3+i_3},\dots,x_{(i_{2_0}+n_2)n_3+i_3}]$}
\STATE{$z:=\text{uFFT(2)}_{r,n_2}(z^T)^T$}
\STATE{$[x_{(i_{2_0}+1)n_3+i_3},\dots,x_{(i_{2_0}+n_2)n_3+i_3}]:=z$}
\ENDFOR
\ENDFOR
\FOR{$i_2=1,\dots,n_2$}
\FOR{$i_3=1,\dots,n_3$}
\STATE{$i_{1_0}:=-1$}
\STATE{$z:=[x_{((i_{1_0}+1)n_2+i_{2}-1)n_3+i_3},\dots,x_{((i_{1_0}+n_1)n_2+i_{2}-1)n_3+i_3}]$}
\STATE{$z:=\text{uFFT(2)}_{r,n_1}(z^T)^T$}
\STATE{$[x_{((i_{1_0}+1)n_2+i_{2}-1)n_3+i_3},\dots,x_{((i_{1_0}+n_1)n_2+i_{2}-1)n_3+i_3}]:=z$}
\ENDFOR
\ENDFOR
\RETURN{$x$}
\end{algorithmic}
\end{algorithm}

\begin{algorithm}[ht]
\caption{Forward FFT, $\text{FFT(2)}_{r,n_1\times\cdots\times n_d}(x)$}
\label{alg:forward-fft}
\begin{algorithmic}[1]
\REQUIRE{$x:=\text{vec}(\mathscr{X})\in\mathbb{C}^{n_1\cdots n_d},\mathscr{X}\in\mathbb{C}^{n_1\times\dots\times n_d}, r\text{ s.t. }\exists\; t_1,\dots,t_d\text{ s.t. }n_1=r^{t_1},\dots,n_d=r^{t_d}$}
\ENSURE{$\text{vec}(\text{DFT}_{n_1\times\dots\times n_d}(\mathscr{X}))$}
\IF{$\text{FFT}_{r,n_1\times\cdots\times n_d}$}
\STATE{$x:=\text{permute}_{r,n_1\times\dots\times n_d}(x)$}
\RETURN{$\text{uFFT}_{r,n_1\times\dots\times n_d}(x)$}
\ELSIF{$\text{FFT2}_{r,n_1\times\cdots\times n_d}$}
\STATE{$x:=\text{uFFT2}_{r,n_1\times\dots\times n_d}(x)$}
\RETURN{$\text{permute}_{r,n_1\times\dots\times n_d}(x)$}
\ENDIF
\end{algorithmic}
\end{algorithm}

\begin{algorithm}[ht]
\caption{Unordered 1D backward FFT, $\text{uiFFT}_{r,n}(y)$}
\label{alg:uiFFT1D}
\begin{algorithmic}[1]
\REQUIRE{$y\in\mathbb{C}^{n}, r\text{ s.t. }\exists\; t\text{ s.t. }n=r^{t}$}
\RETURN{$\overline{A_{r,n}}y/n$}
\COMMENT{See Algos.~\ref{alg:conjugate-butterfly-kernel-radix-2}, \ref{alg:conjugate-butterfly-kernel-radix-4} and \ref{alg:conjugate-butterfly-kernel-radix-8}}
\end{algorithmic}
\end{algorithm}

\begin{algorithm}[ht]
\caption{Unordered 2D backward FFT, $\text{uiFFT}_{r,n_1\times n_2}(y)$}
\label{alg:uiFFT2D}
\begin{algorithmic}[1]
\REQUIRE{$y:=\text{vec}(Y)\in\mathbb{C}^{n_1n_2},\;Y\in\mathbb{C}^{n_1\times n_2},\;r\text{ s.t. }\exists\; t_1,t_2\;\text{ s.t. }\;n_1=r^{t_1},n_2=r^{t_2}$}
\ENSURE{$(\overline{A_{r,n_1}}\otimes\overline{A_{r,n_2}})\, y/(n_1n_2)$}
\FOR{$i_1=1,\dots,n_1$}
\STATE{$i_{2_0}:=(i_1-1)n_2$}
\STATE{$[y_{i_{2_0}+1},\dots,y_{i_{2_0}+n_2}]:=
\text{uiFFT}_{r,n_2}([y_{i_{2_0}+1},\dots,y_{i_{2_0}+n_2}]^T)^T$}
\ENDFOR
\FOR{$i_2=1,\dots,n_2$}
\STATE{$i_{1_0}:=-1$}
\STATE{$z:=[y_{(i_{1_0}+1)n_2+i_2},\dots,y_{(i_{1_0}+n_1)n_2+i_2}]$}
\STATE{$z:=\text{uiFFT}_{r,n_1}(z^T)^T$}
\STATE{$[y_{(i_{1_0}+1)n_2+i_2},\dots,y_{(i_{1_0}+n_1)n_2+i_2}]:=z$}
\ENDFOR
\RETURN{$y/(n_1n_2)$}
\end{algorithmic}
\end{algorithm}

\begin{algorithm}[ht]
\caption{Unordered 3D backward FFT, $\text{uiFFT}_{r,n_1\times n_2\times n_3}(y)$}
\label{alg:uiFFT3D}
\begin{algorithmic}[1]
\REQUIRE{$y:=\text{vec}(\mathscr{Y})\in\mathbb{C}^{n_1n_2n_3},\;\mathscr{Y}\in\mathbb{C}^{n_1\times n_2\times n_3},\;r\text{ s.t. }\exists\; t_1,t_2,t_3\;\text{ s.t. }\;n_1=r^{t_1},n_2=r^{t_2},n_3=r^{t_3}$}
\ENSURE{$(\overline{A_{r,n_1}}\otimes\overline{A_{r,n_2}}\otimes \overline{A_{r,n_3}})\, y/(n_1n_2n_3)$}
\FOR{$i_1=1,\dots,n_1$}
\FOR{$i_2=1,\dots,n_2$}
\STATE{$i_{3_0}:=((i_1-1)n_2+i_2-1)n_3$}
\STATE{$[y_{i_{3_0}+1},\dots,y_{i_{3_0}+n_3}]:=
\text{uiFFT}_{r,n_3}([y_{i_{3_0}+1},\dots,y_{i_{3_0}+n_3}]^T)^T$}
\ENDFOR
\ENDFOR
\FOR{$i_1=1,\dots,n_1$}
\FOR{$i_3=1,\dots,n_3$}
\STATE{$i_{2_0}:=(i_1-1)n_2-1$}
\STATE{$z:=[y_{(i_{2_0}+1)n_3+i_3},\dots,y_{(i_{2_0}+n_2)n_3+i_3}]$}
\STATE{$z:=\text{uiFFT}_{r,n_2}(z^T)^T$}
\STATE{$[y_{(i_{2_0}+1)n_3+i_3},\dots,y_{(i_{2_0}+n_2)n_3+i_3}]:=z$}
\ENDFOR
\ENDFOR
\FOR{$i_2=1,\dots,n_2$}
\FOR{$i_3=1,\dots,n_3$}
\STATE{$i_{1_0}:=-1$}
\STATE{$z:=[y_{((i_{1_0}+1)n_2+i_{2}-1)n_3+i_3},\dots,y_{((i_{1_0}+n_1)n_2+i_{2}-1)n_3+i_3}]$}
\STATE{$z:=\text{uiFFT}_{r,n_1}(z^T)^T$}
\STATE{$[y_{((i_{1_0}+1)n_2+i_{2}-1)n_3+i_3},\dots,y_{((i_{1_0}+n_1)n_2+i_{2}-1)n_3+i_3}]:=z$}
\ENDFOR
\ENDFOR
\RETURN{$y/(n_1n_2n_3)$}
\end{algorithmic}
\end{algorithm}

\begin{algorithm}[ht]
\caption{Backward FFT, $\text{iFFT}_{r,n_1\times\dots\times n_d}(y)$}
\label{alg:backward-fft}
\begin{algorithmic}[1]
\REQUIRE{$y:=\text{vec}(\mathscr{Y})\in\mathbb{C}^{n_1\cdots n_d},\mathscr{Y}\in\mathbb{C}^{n_1\times\dots\times n_d}, r\text{ s.t. }\exists\; t_1,\dots,t_d\text{ s.t. }n_1=r^{t_1},\dots,n_d=r^{t_d}$}
\ENSURE{$\text{vec}(\text{DFT}^{-1}_{n_1\times\dots\times n_d}(\mathscr{Y}))$}
\STATE{$y:=\text{permute}_{r,n_1\times\dots\times n_d}(y)$}
\RETURN{$\text{uiFFT}_{r,n_1\times\dots\times n_d}(y)$}
\end{algorithmic}
\end{algorithm}

\begin{algorithm}[ht]
\caption{Standard FFT-based convolution}
\label{alg:fft-based-convolution}
\begin{algorithmic}[1]
\REQUIRE{$x:=\text{vec}(\mathscr{X})\in\mathbb{C}^{n_1\cdots n_d},\mathscr{X}\in\mathbb{C}^{n_1\times\dots\times n_d},r\text{ s.t. }\exists\; t_1,\dots,t_d\text{ s.t. }n_1=r^{t_1},\dots,n_d=r^{t_d}$}
\ENSURE{$\text{vec}(\mathcal{F}_{\mathscr{H}}(\mathscr{X}))$}
\COMMENT{$\mathscr{H}\in\mathbb{C}^{n_1\times\dots\times n_d}$}
\STATE{$h:=\text{FFT}_{r,n_1\times\dots\times n_d}(h)$}
\COMMENT{$h:=\text{vec}(\mathscr{H})\in\mathbb{C}^{n_1\cdots n_d}$}
\STATE{$x:=\text{FFT}_{r,n_1\times\dots\times n_d}(x)$}
\STATE{$x:=h\circ x$}
\RETURN{$\text{iFFT}_{r,n_1\times\dots\times n_d}(x)$}
\end{algorithmic}
\end{algorithm}

\begin{algorithm}[ht]
\caption{Permutation-avoiding convolution}
\label{alg:permutation-avoiding-convolution}
\begin{algorithmic}[1]
\REQUIRE{$x:=\text{vec}(\mathscr{X})\in\mathbb{C}^{n_1\dots n_d},\mathscr{X}\in\mathbb{C}^{n_1\times\dots\times n_d}, r\text{ s.t. }\exists\; t_1,\dots,t_d\text{ s.t. }n_1=r^{t_1},\dots,n_d=r^{t_d}$}
\ENSURE{$\text{vec}(\mathcal{F}_{\mathscr{H}}(\mathscr{X}))$}
\COMMENT{$\mathscr{H}\in\mathbb{C}^{n_1\times\dots\times n_d}$}
\STATE{$h:=\text{uFFT2}_{r,n_1\times\dots\times n_d}(h)$}
\COMMENT{$h:=\text{vec}(\mathscr{H})\in\mathbb{C}^{n_1\cdots n_d}$}
\STATE{$x:=\text{uFFT2}_{r,n_1\times\dots\times n_d}(x)$}
\STATE{$x:=h\circ x$}
\RETURN{$\text{uiFFT}_{r,n_1\times\dots\times n_d}(x)$}
\end{algorithmic}
\end{algorithm}

To characterize the performance regime of the algorithms presented within the roofline model framework, we assess theoretical estimates of arithmetic intensity.
The FLOP count for FFT operations stems entirely from butterfly operations, yielding $5n\log_2n$, $4.25n\log_2n$, and $4.08n\log_2n$ FLOP for 1D radix-2, radix-4, and radix-8 implementations, respectively.
While simplistic and incomplete, counting memory read and write accesses provides an initial estimate of memory traffic patterns.
In 1D implementations, each of the $\log_r n$ butterfly stages reads and writes $n$ elements.
The permutation operation contributes an additional $n$ reads and $n$ writes.
This totals $2n\cdot(1+\log_r n)$ complex double-precision elements.
For each FFT kernel, the theoretical arithmetic intensity becomes:
\begin{itemize}
\item Radix-2: 
$(5n \log_2 n~\text{FLOP})/(2n\cdot(1+\log_2 n) \times 16~\text{byte}) 
=
0.156(\log_2n)/(1+\log_2 n)~\text{FLOP/byte}$,
\item Radix-4: 
$(4.25n \log_2 n~\text{FLOP})/(2n\cdot(1+\log_4 n) \times 16~\text{byte}) 
= 
0.133(\log_2 n)/(1+\log_4 n)~\text{FLOP/byte}$,
\item Radix-8: 
$(4.08n \log_2 n~\text{FLOP})/(2n\cdot(1+\log_8 n) \times 16~\text{byte}) 
= 
0.128(\log_2 n)/(1+\log_8 n)~\text{FLOP/byte}$.
\end{itemize}
These expressions indicate arithmetic intensity that increases with problem size $n$ towards limits of $0.156$, $0.266$ and $0.383$~FLOP/byte, respectively.
Permutation-avoiding implementations reach these limits by reducing memory traffic to $2n\log_r n$ accesses.
Indeed, we have:
\begin{itemize}
\item Radix-2: $(5n \log_2 n~\text{FLOP})/(2n\log_2 n \times 16~\text{byte}) = 0.156~\text{FLOP/byte}$,
\item Radix-4: $(4.25n \log_2 n~\text{FLOP})/(2n\log_4 n \times 16~\text{byte}) = 0.266~\text{FLOP/byte}$,
\item Radix-8: $(4.08n \log_2 n~\text{FLOP})/(2n\log_8 n \times 16~\text{byte}) = 0.383~\text{FLOP/byte}$.
\end{itemize}
Theoretically, permutation-avoiding implementations maintain constant arithmetic intensity independent of transform size.
However, this analysis oversimplifies the memory hierarchy's complexity.
In practice, data movement occurs across multiple memory levels, with performance determined by the algorithm's temporal and spatial locality characteristics and the underlying memory architecture.
As transform size $n$ increases, locality degrades and memory traffic intensifies when data exceeds the capacity of successive cache levels.
Consequently, profiling DRAM traffic provides a more representative measure of actual data movement costs.

\section{Numerical experiments}\label{sec:numerical-experiments}
In our experiments, we benchmark and verify our implementations against FFTW using the \texttt{FFTW\_ESTIMATE} flag.
Since FFTW under \texttt{FFTW\_ESTIMATE} does not use empirical performance tuning, this setting offers the most inherently fair comparison with our implementation.
By benchmarking both the complete FFT pipeline and isolated operations (index-reversal permutation and butterfly computations), we obtain fine-grained insights that highlight how a clean, portable general radix implementation performs relative to FFTW's general-purpose performance without platform-specific optimizations.

All experiments are conducted on two distinct computing systems:
\paragraph{Computing system \#1} Acer Swift SFL14-53M laptop with the following specifications:
\begin{itemize}
\item Processor: Intel Core i5-1334U (13th Generation Raptor Lake)
\begin{itemize}
\item Architecture: x86\_64 hybrid (2 performance cores + 8 efficiency cores, 12 threads total)  
\item Base/Boost frequency: 400 MHz - 4.6 GHz (maximum)
\item Cache hierarchy: L1d 352 KiB, L1i 576 KiB, L2 6.5 MiB, L3 12 MiB
\item Instruction sets: AVX, AVX2, AVX\_VNNI, SSE 4.2
\end{itemize}
\item Memory: $2\times32$ GiB DDR5-4800, peak bandwidth measured at 14.55~GB/s per core.
\item Software Environment:
\begin{itemize}
\item Operating System: Ubuntu 24.04.4 LTS (Noble Numbat)
\item Kernel: Linux 6.17.0-20-generic  
\item Compiler: GCC 13.3.0
\end{itemize}
\end{itemize}

\paragraph{Computing system \#2} Dual-socket Intel Xeon server with the following specifications:
\begin{itemize}
\item Processor: Intel Xeon Platinum 8480+ (4th Generation Sapphire Rapids)
\begin{itemize}
\item Architecture: x86\_64 dual-socket (2 × 56 cores, 224 threads total)
\item Base/Boost frequency: 800 MHz - 3.8 GHz
\item Cache hierarchy: L1 8.8 MiB, L2 224 MiB, L3 210 MiB
\item Instruction sets: AVX, AVX2, AVX-512F/DQ/BW/VL, AMX (Advanced Matrix Extensions)
\end{itemize}
\item Memory: 503 GiB DDR5, peak bandwidth measured at 9.17~GB/s per core.
\item Software Environment:
\begin{itemize}
\item Operating System: Red Hat Enterprise Linux 9.5
\item Kernel: Linux 5.14.0-503.40.1.el9\_5.x86\_64
\item Compiler: GCC 11.5.0
\end{itemize}
\end{itemize}

The 12 MiB L3 cache of the first computing system accommodates moderate-sized FFT transforms entirely in cache, whereas the massive 210 MiB L3 cache of the second system can accommodate very large transforms entirely in cache.

With these technical specifications in hand, one can assess machine balance for each system and compare these values with the theoretical arithmetic intensities presented in Section~\ref{sec:algorithms} to evaluate whether, in theory, each kernel is compute-bound or memory-bound.
In particular, assuming single-thread performance without explicit vectorization or specific instruction set optimizations, computing system~\#1 achieves a peak performance of 4.6~GFLOP/s and a single-core peak memory bandwidth of 14.55~GB/s.
This yields a machine balance of 0.316~FLOP/byte, which suggests that radix-2 and radix-4 implementations are memory-bound, both with and without permutations. 
However, permutation-avoiding kernels of radix-8 operate in the compute-bound regime, irrespective of $n$, whereas full-FFT radix-8 kernels become compute-bound for $n\geq 8^{5}$.
For computing system~\#2, we observe a single-core peak performance of 3.8~GFLOP/s and a peak memory bandwidth of 9.17~GB/s per core, resulting in a machine balance of 0.414~FLOP/byte.
Under these conditions, all kernel implementations are memory-bound.

\subsection{1D results}
Tables~\ref{tab:1d-fft} and \ref{tab:1d-ifft} present forward and backward FFT timing results on system~\#1, respectively.
Runtimes measured on system~\#2 are presented in Tables~\ref{tab:1d-fft-2} and \ref{tab:1d-ifft-2}.
Several key observations emerge from these results.

Despite radix-8 requiring the fewest FLOP (4.08 $\log_2 n$ vs 4.25 $\log_2 n$ for radix-4 and 5 $\log_2 n$ for radix-2), radix-4 consistently delivers the best performance across all problem sizes, on both systems. 
This counterintuitive result demonstrates that theoretical FLOP counts provide an incomplete picture of real-world performance.
The poor performance of radix-8 can be attributed to scattered memory access patterns across 8 elements with large strides, and increased register pressure from numerous temporary variables. 
These implementation overheads overwhelm the theoretical FLOP advantage.
For radix-4 implementations, index-reversal permutations account for 48--71\% of total transform time on system \#1, and 55--73\% on system \#2, indicating significant potential speedup for permutation-avoiding implementations. 
On system \#1, our radix-4 implementation of forward and backward FFT kernels achieves runtime ranging from $0.71\times$ to $1.44\times$ that of FFTW, generally doing better than FFTW for larger transform sizes, whereas FFTW outperforms our implementation for smaller problem sizes.
On system \#2, FFTW consistently outperforms our radix-4 implementations with speedups ranging from $1.34\times$ to $2.36\times$.
Overall, the performance of our FFT implementations is on par with that of state-of-the-art libraries when those do not leverage architecture-specific optimizations and empirical transform planning.
This observation validates the findings of our numerical experiments performed to investigate the performance of permutation-avoiding convolution kernels.

Tables~\ref{tab:1d-convolution} and \ref{tab:1d-convolution-2} present the runtime of 1D FFT-based and permutation-avoiding convolution kernels on system \#1 and \#2, respectively.
Similarly to standard FFTs, our best-performing FFT-based and permutation-avoiding convolution kernels are those with radix-4.
On system \#1, the switch from a standard radix-4 FFT-based convolution to a permutation-avoiding kernel leads to consistent speedups ranging from $2.30\times$ to $4.04\times$.
On system \#2, those numbers range from $2.20\times$ to $3.97\times$.

\begin{table}
\caption{Timing of 1D index-reversal permutation, unordered and full forward FFT --- Computing system \#1}
\label{tab:1d-fft}
\begin{tabular}{cccccc}
\toprule
Procedure     & $n=2^{20}$ & $n=2^{22}$ & $n=2^{24}$ & 
$n=2^{26}$ & $n=2^{28}$ \\
\midrule
\multicolumn{6}{c}{radix-2}\\
\shortstack{Index-reversal\\permutation} 
              &  0.011~s & 0.067~s & 0.270~s & 1.147~s & 10.433~s\\
Unordered FFT &  0.015~s & 0.075~s & 0.334~s & 1.436~s &  6.770~s\\
FFT           &  0.027~s & 0.142~s & 0.602~s & 2.582~s & 17.205~s\\
\midrule
\multicolumn{6}{c}{radix-4}\\
\shortstack{Index-reversal\\permutation}   
              &  0.011~s & 0.064~s & 0.271~s & 1.149~s & 10.445~s\\
Unordered FFT &  0.011~s & 0.050~s & 0.219~s & 0.957~s &  4.336~s\\
FFT           &  0.023~s & 0.115~s & 0.492~s & 2.127~s & 14.777~s\\
\midrule
\multicolumn{6}{c}{radix-8}\\
\shortstack{Index-reversal\\permutation}
              &  -       & -       & 0.270~s &      -  & -      \\
Unordered FFT &  -       & -       & 0.652~s &      -  & -      \\
FFT           &  -       & -       & 0.922~s &      -  & -      \\
\midrule
FFT (FFTW)    &  0.016~s & 0.080~s & 0.595~s & 2.752~s & 12.539~s\\
\bottomrule
\end{tabular}
\end{table}

\begin{table}
\caption{Timing of 1D index-reversal permutation, unordered and full backward FFT --- Computing system \#1}
\label{tab:1d-ifft}
\begin{tabular}{cccccc}
\toprule
Procedure     & $n=2^{20}$ & $n=2^{22}$ & $n=2^{24}$ & 
$n=2^{26}$ & $n=2^{28}$ \\
\midrule
\multicolumn{6}{c}{radix-2}\\
\shortstack{Index-reversal\\permutation}
              &  0.011~s & 0.065~s & 0.269~s & 1.151~s &  5.115~s\\
Unordered FFT &  0.016~s & 0.077~s & 0.338~s & 1.474~s &  6.989~s\\
FFT           &  0.027~s & 0.142~s & 0.607~s & 2.626~s & 12.106~s\\
\midrule
\multicolumn{6}{c}{radix-4}\\
\shortstack{Index-reversal\\permutation}
              &  0.011~s & 0.064~s & 0.275~s & 1.156~s &  5.141~s\\
Unordered FFT &  0.012~s & 0.053~s & 0.256~s & 1.042~s &  4.525~s\\
FFT           &  0.023~s & 0.117~s & 0.526~s & 2.185~s &  9.673~s\\
\midrule
\multicolumn{6}{c}{radix-8}\\
\shortstack{Index-reversal\\permutation}
              &  -       &   -     & 0.267~s &    -    & -\\
Unordered FFT &  -       &   -     & 0.681~s &    -    & -\\
FFT           &  -       &   -     & 0.948~s &    -    & -\\
\midrule
FFT (FFTW)    &  0.016~s & 0.083~s & 0.693~s & 3.014~s & 12.723~s\\
\bottomrule
\end{tabular}
\end{table}

\begin{table}
\caption{Timing of 1D index-reversal permutation, unordered and full forward FFT --- Computing system \#2}
\label{tab:1d-fft-2}
\begin{tabular}{cccccc}
\toprule
Procedure     & $n=2^{20}$ & $n=2^{22}$ & $n=2^{24}$ & 
$n=2^{26}$ & $n=2^{28}$ \\
\midrule
\multicolumn{6}{c}{radix-2}\\
\shortstack{Index-reversal\\permutation} 
              &  0.017~s & 0.102~s & 0.860~s & 4.167~s & 17.507~s\\
Unordered FFT &  0.014~s & 0.068~s & 0.499~s & 2.477~s & 11.086~s\\
FFT           &  0.030~s & 0.172~s & 1.354~s & 6.645~s & 28.857~s\\
\midrule
\multicolumn{6}{c}{radix-4}\\
\shortstack{Index-reversal\\permutation}   
              &  0.014~s & 0.096~s & 0.848~s & 4.828~s & 18.015~s\\
Unordered FFT &  0.011~s & 0.053~s & 0.314~s & 1.752~s &  7.118~s\\
FFT           &  0.025~s & 0.152~s & 1.157~s & 6.637~s & 25.332~s\\
\midrule
\multicolumn{6}{c}{radix-8}\\
\shortstack{Index-reversal\\permutation}
              &  -       & -       & 0.866~s &      -  & -      \\
Unordered FFT &  -       & -       & 0.694~s &      -  & -      \\
FFT           &  -       & -       & 1.556~s &      -  & -      \\
\midrule
FFT (FFTW)    &  0.016~s & 0.100~s & 0.669~s & 2.815~s & 11.915~s\\
\bottomrule
\end{tabular}
\end{table}

\begin{table}
\caption{Timing of 1D index-reversal permutation, unordered and full backward FFT --- Computing system \#2}
\label{tab:1d-ifft-2}
\begin{tabular}{cccccc}
\toprule
Procedure     & $n=2^{20}$ & $n=2^{22}$ & $n=2^{24}$ & 
$n=2^{26}$ & $n=2^{28}$ \\
\midrule
\multicolumn{6}{c}{radix-2}\\
\shortstack{Index-reversal\\permutation}
              &  0.018~s & 0.093~s & 1.009~s & 4.402~s & 17.793~s\\
Unordered FFT &  0.015~s & 0.070~s & 0.573~s & 2.649~s & 11.680~s\\
FFT           &  0.033~s & 0.173~s & 1.594~s & 7.067~s & 29.472~s\\
\midrule
\multicolumn{6}{c}{radix-4}\\
\shortstack{Index-reversal\\permutation}
              &  0.015~s & 0.096~s & 0.959~s & 4.356~s & 17.824~s\\
Unordered FFT &  0.012~s & 0.057~s & 0.356~s & 1.641~s &  7.305~s\\
FFT           &  0.027~s & 0.176~s & 1.328~s & 6.005~s & 25.123~s\\
\midrule
\multicolumn{6}{c}{radix-8}\\
\shortstack{Index-reversal\\permutation}
              &  -       &   -     & 0.990~s &    -    & -\\
Unordered FFT &  -       &   -     & 0.724~s &    -    & -\\
FFT           &  -       &   -     & 1.715~s &    -    & -\\
\midrule
FFT (FFTW)    &  0.017~s & 0.131~s & 0.720~s & 2.919~s & 12.513~s\\
\bottomrule
\end{tabular}
\end{table}

\begin{table}
\caption{Timing of 1D permutation-avoiding and standard FFT-based convolution --- Computing system \#1}
\label{tab:1d-convolution}
\begin{tabular}{cccccc}
\toprule
Procedure     & $n=2^{20}$ & $n=2^{22}$ & $n=2^{24}$ & 
$n=2^{26}$ & $n=2^{28}$ \\
\midrule
\multicolumn{6}{c}{radix-2}\\
\shortstack{Permutation-avoiding\\convolution$\vphantom{x_{|_{|_|}}}$} &  
0.046~s & 0.226~s & 0.997~s & 4.323~s & 25.051~s\\
\shortstack{Standard FFT-based\\convolution} &  
0.084~s & 0.431~s & 1.845~s & 7.991~s & 62.357~s\\
\midrule
\multicolumn{6}{c}{radix-4}\\
\shortstack{Permutation-avoiding\\convolution$\vphantom{x_{|_{|_|}}}$} &  
0.033~s & 0.149~s & 0.654~s & 2.841~s & 13.603~s\\
\shortstack{Standard FFT-based\\convolution} &  
0.076~s & 0.352~s & 1.508~s & 6.570~s & 54.874~s\\
\midrule
\multicolumn{6}{c}{radix-8}\\
\shortstack{Permutation-avoiding\\convolution$\vphantom{x_{|_{|_|}}}$} &  
   -    &    -    & 2.027~s &     -    & -\\
\shortstack{Standard FFT-based\\convolution} &  
   -    &    -    & 2.844~s &     -    & -\\
\midrule
\shortstack{Standard FFT-based\\convolution (FFTW)} &
0.127~s & 0.257~s & 1.872~s & 8.642~s & 38.620~s\\
\bottomrule
\end{tabular}
\end{table}

\begin{table}
\caption{Timing of 1D permutation-avoiding and standard FFT-based convolution --- Computing system \#2}
\label{tab:1d-convolution-2}
\begin{tabular}{cccccc}
\toprule
Procedure     & $n=2^{20}$ & $n=2^{22}$ & $n=2^{24}$ & 
$n=2^{26}$ & $n=2^{28}$ \\
\midrule
\multicolumn{6}{c}{radix-2}\\
\shortstack{Permutation-avoiding\\convolution$\vphantom{x_{|_{|_|}}}$} &  
0.049~s & 0.294~s & 1.964~s &  9.009~s & 39.872~s\\
\shortstack{Standard FFT-based\\convolution} &  
0.093~s & 0.729~s & 4.906~s & 21.483~s & 89.465~s\\
\midrule
\multicolumn{6}{c}{radix-4}\\
\shortstack{Permutation-avoiding\\convolution$\vphantom{x_{|_{|_|}}}$} &  
0.035~s & 0.186~s & 1.052~s &  4.752~s & 21.036~s\\
\shortstack{Standard FFT-based\\convolution} &  
0.077~s & 0.665~s & 4.180~s & 18.337~s & 76.572~s\\
\midrule
\multicolumn{6}{c}{radix-8}\\
\shortstack{Permutation-avoiding\\convolution$\vphantom{x_{|_{|_|}}}$} &  
   -    &    -    & 2.103~s &     -    & -\\
\shortstack{Standard FFT-based\\convolution} &  
   -    &    -    & 5.249~s &     -    & -\\
\midrule
\shortstack{Standard FFT-based\\convolution (FFTW)} &
0.056~s & 0.375~s & 2.246~s & 9.258~s & 39.281~s\\
\bottomrule
\end{tabular}
\end{table}

\subsection{2D results}
Tables~\ref{tab:2d-fft}, \ref{tab:2d-ifft}, \ref{tab:2d-fft-2} and \ref{tab:2d-ifft-2} show isotropic 2D FFT results, that is, of transforms with $n_1=n_2\in\{10,11,12,13,14\}$.
Note that those total transform sizes $n=n_1n_2$ equal those of the 1D transforms we presented in the last Section, that is, $n=2^{t}$ for $t\in\{20,22,24,26,28\}$. 
We observe that the application of precomputed 2D index-reversal permutations of the form $P_{r,n_1}\otimes P_{r,n_2}$ achieves much higher data locality than equally sized 1D permutations $P_{r,n_1n_2}$, resulting in consistent speedups of the index reversal operation spanning from 2.72$\times$ to 9.14$\times$ on system \#1, and 7.00$\times$ to 14.53$\times$ on system \#2.
For the total transform sizes mentioned above, $n_1$ and $n_2$ are more indicative of the stride of 2D permutations than $n=n_1n_2$.
Indeed, arrays of size $\sqrt{n}$ easily fit in L1 cache, on both computing systems, which explains the observed speedups of isotropic transforms.
In an attempt to better illustrate the effect of data traffic induced by 2D index-reversal permutation on runtime, in addition to FFT-based and permutation-avoiding convolution kernels of isotropic transform sizes (see Tables~\ref{tab:2d-convolution} and \ref{tab:2d-convolution-2} for systems \#1 and \#2, respectively), we also present results for anisotropic transform sizes, i.e., with $n_1\gg n_2$ for which arrays of size $n_1$ do not fit in cache (see Tables~\ref{tab:2d-convolution-anisotropic} and \ref{tab:2d-convolution-anisotropic-2} in appendix).
But even then, the relative contribution of 2D butterfly operations to runtime is much more significant than in 1D.
This is explained by the fact that 2D butterfly operations also induce non-negligible data traffic due to the gathering and scattering of column data, thus reducing effective cache utilization during the separable transform process.
As a result, the speedup achieved by permutation-avoiding over FFT-based convolution kernels are 1.12--1.37$\times$ on system \#1, and 1.02--1.27$\times$ on system \#2. 
As reported in those tables, our implementations remain competitive with FFTW, often matching or exceeding FFTW's performance for larger problem sizes, particularly for radix-4.

\begin{table}
\caption{Timing of isotropic 2D permutation-avoiding and standard FFT-based convolution --- Computing system \#1}
\label{tab:2d-convolution}
\begin{tabular}{cccccc}
\toprule
Procedure     & 
$n_1=n_2=2^{10}$& 
$n_1=n_2=2^{11}$& 
$n_1=n_2=2^{12}$& 
$n_1=n_2=2^{13}$& 
$n_1=n_2=2^{14}$\\
\midrule
\multicolumn{6}{c}{radix-2}\\
\shortstack{Permutation-avoiding\\convolution$\vphantom{x_{|_{|_|}}}$} &  
0.084~s & 0.589~s & 3.158~s & 13.551~s & 59.071~s\\
\shortstack{Standard FFT-based\\convolution} &  
0.085~s & 0.607~s & 3.405~s & 14.704~s & 64.071~s\\
\midrule
\multicolumn{6}{c}{radix-4}\\
\shortstack{Permutation-avoiding\\convolution$\vphantom{x_{|_{|_|}}}$} &  
0.087~s & -       & 3.130~s & -        & 57.637~s\\
\shortstack{Standard FFT-based\\convolution} &  
0.090~s & -       & 3.446~s & -        & 63.018~s\\
\midrule
\shortstack{Standard FFT-based\\convolution (FFTW)} &
0.065~s & 0.414~s & 2.770~s & 13.930~s & 100.659~s\\
\bottomrule
\end{tabular}
\end{table}

\begin{table}
\caption{Timing of isotropic 2D permutation-avoiding and standard FFT-based convolution --- Computing system \#2}
\label{tab:2d-convolution-2}
\begin{tabular}{cccccc}
\toprule
Procedure     & 
$n_1=n_2=2^{10}$& 
$n_1=n_2=2^{11}$& 
$n_1=n_2=2^{12}$& 
$n_1=n_2=2^{13}$& 
$n_1=n_2=2^{14}$\\
\midrule
\multicolumn{6}{c}{radix-2}\\
\shortstack{Permutation-avoiding\\convolution$\vphantom{x_{|_{|_|}}}$} &  
0.106~s & 0.444s & 3.436~s &  17.812~s & 73.649~s\\
\shortstack{Standard FFT-based\\convolution} &  
0.113~s & 0.461s & 3.667~s &  18.819~s & 78.202~s\\
\midrule
\multicolumn{6}{c}{radix-4}\\
\shortstack{Permutation-avoiding\\convolution$\vphantom{x_{|_{|_|}}}$} &  
0.104~s & -       & 3.371~s & -        & 73.019~s\\
\shortstack{Standard FFT-based\\convolution} &  
0.108~s & -       & 3.571~s & -        & 77.105~s\\
\midrule
\shortstack{Standard FFT-based\\convolution (FFTW)} &
0.099~s & 0.555~s & 4.610~s & 22.525~s & 127.934~s\\
\bottomrule
\end{tabular}
\end{table}

\begin{table}
\caption{Timing of anisotropic 2D permutation-avoiding and standard FFT-based convolution --- Computing system \#1}
\label{tab:2d-convolution-anisotropic}
\begin{tabular}{ccccc}
\toprule
Procedure     & 
\shortstack{$n_1=2^{20}$\\$n_2=2$}& 
\shortstack{$n_1=2^{22}$\\$n_2=2$}& 
\shortstack{$n_1=2^{24}$\\$n_2=2$}& 
\shortstack{$n_1=2^{26}$\\$n_2=2$}\\
\midrule
&\multicolumn{4}{c}{radix-2}\\
\shortstack{Permutation-avoiding\\convolution$\vphantom{x_{|_{|_|}}}$} &  
0.330~s & 1.444~s & 6.162~s & 30.759~s \\
\shortstack{Standard FFT-based\\convolution} &  
0.409~s & 1.771~s & 7.589~s & 34.558~s \\
\midrule
&\multicolumn{4}{c}{radix-4}\\
\shortstack{Permutation-avoiding\\convolution$\vphantom{x_{|_{|_|}}}$} &  
0.235~s & 1.049~s & 4.512~s & 19.841   \\
\shortstack{Standard FFT-based\\convolution} &  
0.322~s & 1.391~s & 5.993~s & 26.155~s \\
\midrule
\shortstack{Standard FFT-based\\convolution (FFTW)} &
0.305~s & 1.494~s & 7.031~s & 34.566~s \\
\bottomrule
\end{tabular}
\end{table}

\begin{table}
\caption{Timing of anisotropic 2D permutation-avoiding and standard FFT-based convolution --- Computing system \#2}
\label{tab:2d-convolution-anisotropic-2}
\begin{tabular}{ccccc}
\toprule
Procedure   & 
\shortstack{$n_1=2^{20}$\\$n_2=2$}& 
\shortstack{$n_1=2^{22}$\\$n_2=2$}& 
\shortstack{$n_1=2^{24}$\\$n_2=2$}& 
\shortstack{$n_1=2^{26}$\\$n_2=2$}\\
\midrule
&\multicolumn{4}{c}{radix-2}\\
\shortstack{Permutation-avoiding\\convolution$\vphantom{x_{|_{|_|}}}$} &  
0.604~s & 1.938~s & 11.099~s & 50.402~s \\
\shortstack{Standard FFT-based\\convolution} &  
0.619~s & 2.133~s & 11.819~s & 52.703~s \\
\midrule
&\multicolumn{4}{c}{radix-4}\\
\shortstack{Permutation-avoiding\\convolution$\vphantom{x_{|_{|_|}}}$} &  
0.339~s & 1.486~s & 7.546~s & 33.361   \\
\shortstack{Standard FFT-based\\convolution} &  
0.386~s & 1.886~s & 9.333~s & 40.698~s \\
\midrule
\shortstack{Standard FFT-based\\convolution (FFTW)} &
0.320~s & 2.572~s & 10.436~s & 46.978~s \\
\bottomrule
\end{tabular}
\end{table}

\subsection{3D results}
Tables~\ref{tab:3d-fft}, \ref{tab:3d-ifft}, \ref{tab:3d-fft-2} and \ref{tab:3d-ifft-2} present isotropic 3D FFT results, i.e., with $n_1=n_2=n_3$.
The phenomena discussed in the last Section about 2D transforms are magnified in 3D.
That is, the data locality of precomputed 3D index-reversal permutations is better characterized by $n_1,n_2$ and $n_3$ than by $n=n_1n_2n_3$, resulting in even smaller fractions of the runtime of each transform, making the achievable speedup by permutation-avoiding implementations less significant than in 2D and even more so compared to the 1D case.
In Tables~\ref{tab:3d-convolution-anisotropic} and \ref{tab:3d-convolution-anisotropic-2}, we report the runtimes of permutation-avoiding and FFT-based convolution on anisotropic data, i.e., with $n_1\gg n_2,n_3$ such that arrays of size $n_1$ do not fit in cache.
Even then, we note that the speedup achieved by permutation-avoidance is less significant than in 2D.
Contributing to this phenomenon is even further data traffic induced by the fetching and scattering of data in one more dimension when apply 1D butterfly operations in each direction.
As reported in those tables, our 3D implementations are competitive with FFTW, which validates our observations on the negligible speedup achievable by permutation-avoidance in higher-dimensional settings.

\begin{table}
\caption{Timing of isotropic 3D permutation-avoiding and standard FFT-based convolution --- Computing system \#1}
\label{tab:3d-convolution}
\begin{tabular}{cccc}
\toprule
Procedure     & 
$n_1=n_2=n_3=2^{7}$& 
$n_1=n_2=n_3=2^{8}$& 
$n_1=n_2=n_3=2^{9}$\\
\midrule
\multicolumn{4}{c}{radix-2}\\
\shortstack{Permutation-avoiding\\convolution$\vphantom{x_{|_{|_|}}}$} &  
0.148~s & 1.709~s & 20.745~s\\
\shortstack{Standard FFT-based\\convolution} &  
0.153~s & 1.752~s & 21.168~s\\
\midrule
\multicolumn{4}{c}{radix-4}\\
\shortstack{Permutation-avoiding\\convolution$\vphantom{x_{|_{|_|}}}$} &  
-       & 1.644~s & -      \\
\shortstack{Standard FFT-based\\convolution} &  
-       & 1.716~s & -      \\
\midrule
\shortstack{Standard FFT-based\\convolution (FFTW)} &
0.053~s & 1.530~s & 23.984~s\\
\bottomrule
\end{tabular}
\end{table}

\begin{table}
\caption{Timing of isotropic 3D permutation-avoiding and standard FFT-based convolution --- Computing system \#2}
\label{tab:3d-convolution-2}
\begin{tabular}{cccc}
\toprule
Procedure     & 
$n_1=n_2=n_3=2^{7}$& 
$n_1=n_2=n_3=2^{8}$& 
$n_1=n_2=n_3=2^{9}$\\
\midrule
\multicolumn{4}{c}{radix-2}\\
\shortstack{Permutation-avoiding\\convolution$\vphantom{x_{|_{|_|}}}$} &  
0.245~s & 2.677~s & 31.650~s\\
\shortstack{Standard FFT-based\\convolution} &  
0.252~s & 2.762~s & 32.303~s\\
\midrule
\multicolumn{4}{c}{radix-4}\\
\shortstack{Permutation-avoiding\\convolution$\vphantom{x_{|_{|_|}}}$} &  
-       & 2.477~s & -      \\
\shortstack{Standard FFT-based\\convolution} &  
-       & 2.530~s & -      \\
\midrule
\shortstack{Standard FFT-based\\convolution (FFTW)} &
0.053~s & 3.146~s & 35.366~s\\
\bottomrule
\end{tabular}
\end{table}

\begin{table}
\caption{Timing of anisotropic 3D permutation-avoiding and standard FFT-based convolution --- Computing system \#1}
\label{tab:3d-convolution-anisotropic}
\begin{tabular}{cccc}
\toprule
Procedure     & 
\shortstack{$n_1=2^{20}$\\$n_2=n_3=2$}& 
\shortstack{$n_1=2^{22}$\\$n_2=n_3=2$}& 
\shortstack{$n_1=2^{24}$\\$n_2=n_3=2$}\\
\midrule
&\multicolumn{3}{c}{radix-2}\\
\shortstack{Permutation-avoiding\\convolution$\vphantom{x_{|_{|_|}}}$} &  
2.423~s & 10.277~s & 44.562~s\\
\shortstack{Standard FFT-based\\convolution} &  
2.658~s & 11.347~s & 48.183~s\\
\midrule
&\multicolumn{3}{c}{radix-4}\\
\shortstack{Permutation-avoiding\\convolution$\vphantom{x_{|_{|_|}}}$} &  
2.067~s & 8.295~s & 34.198~s\\
\shortstack{Standard FFT-based\\convolution} &  
2.334~s & 9.413~s & 39.065~s\\
\midrule
\shortstack{Standard FFT-based\\convolution (FFTW)} &
2.213~s & 10.673~s & 51.115~s\\
\bottomrule
\end{tabular}
\end{table}

\begin{table}
\caption{Timing of anisotropic 3D permutation-avoiding and standard FFT-based convolution --- Computing system \#2}
\label{tab:3d-convolution-anisotropic-2}
\begin{tabular}{cccc}
\toprule
Procedure   &
\shortstack{$n_1=2^{20}$\\$n_2=n_3=2$}&
\shortstack{$n_1=2^{22}$\\$n_2=n_3=2$}& 
\shortstack{$n_1=2^{24}$\\$n_2=n_3=2$}\\
\midrule
&\multicolumn{3}{c}{radix-2}\\
\shortstack{Permutation-avoiding\\convolution$\vphantom{x_{|_{|_|}}}$} &  
4.373~s & 16.097~s & 75.846~s\\
\shortstack{Standard FFT-based\\convolution} &  
4.564~s & 16.107~s & 75.525~s\\
\midrule
&\multicolumn{3}{c}{radix-4}\\
\shortstack{Permutation-avoiding\\convolution$\vphantom{x_{|_{|_|}}}$} &  
3.379~s & 13.917~s & 60.246~s\\
\shortstack{Standard FFT-based\\convolution} &  
3.450~s & 14.805~s & 64.209~s\\
\midrule
\shortstack{Standard FFT-based\\convolution (FFTW)} &
3.759~s & 17.464~s & 77.078~s\\
\bottomrule
\end{tabular}
\end{table}

\subsection{Scaling analysis}
\begin{figure}[h]
\centering
\includegraphics[width=\linewidth]{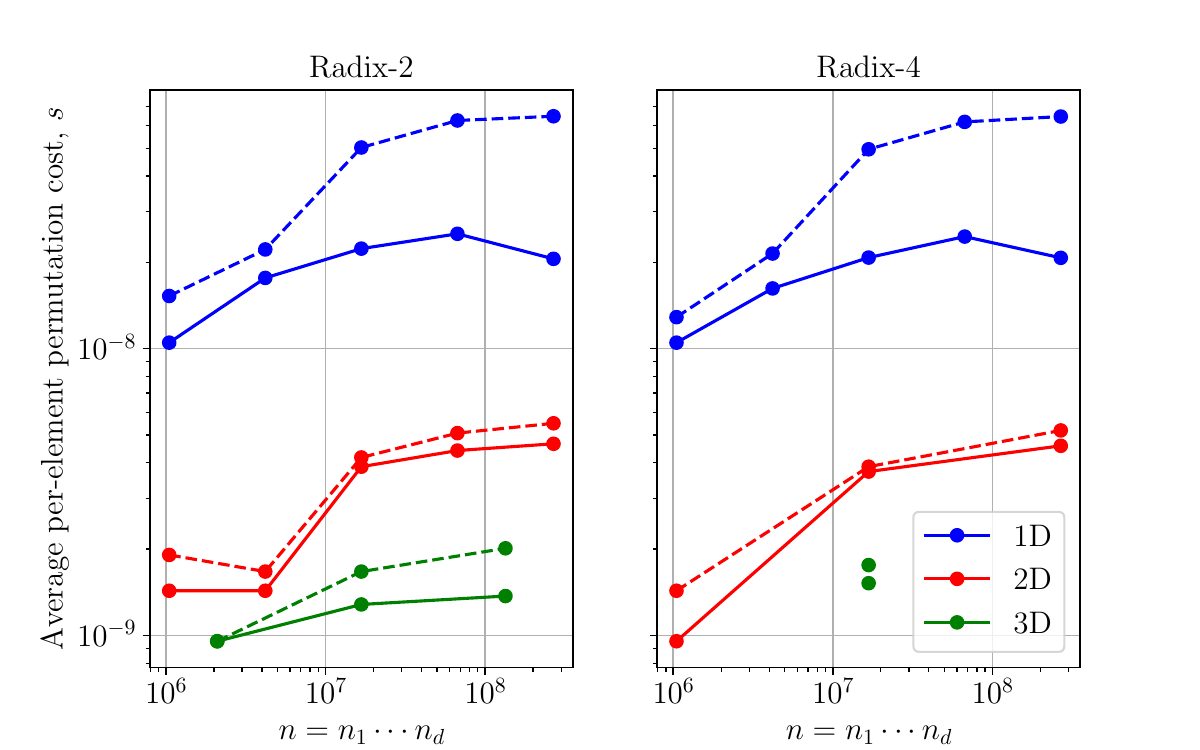}
\caption{Average per-element permutation cost for computing systems \#1 (solid lines) and \#2 (dashed lines)}
\label{fig:per-element-permutation-cost}
\end{figure}

\begin{figure}[h]
\centering
\includegraphics[width=\linewidth]{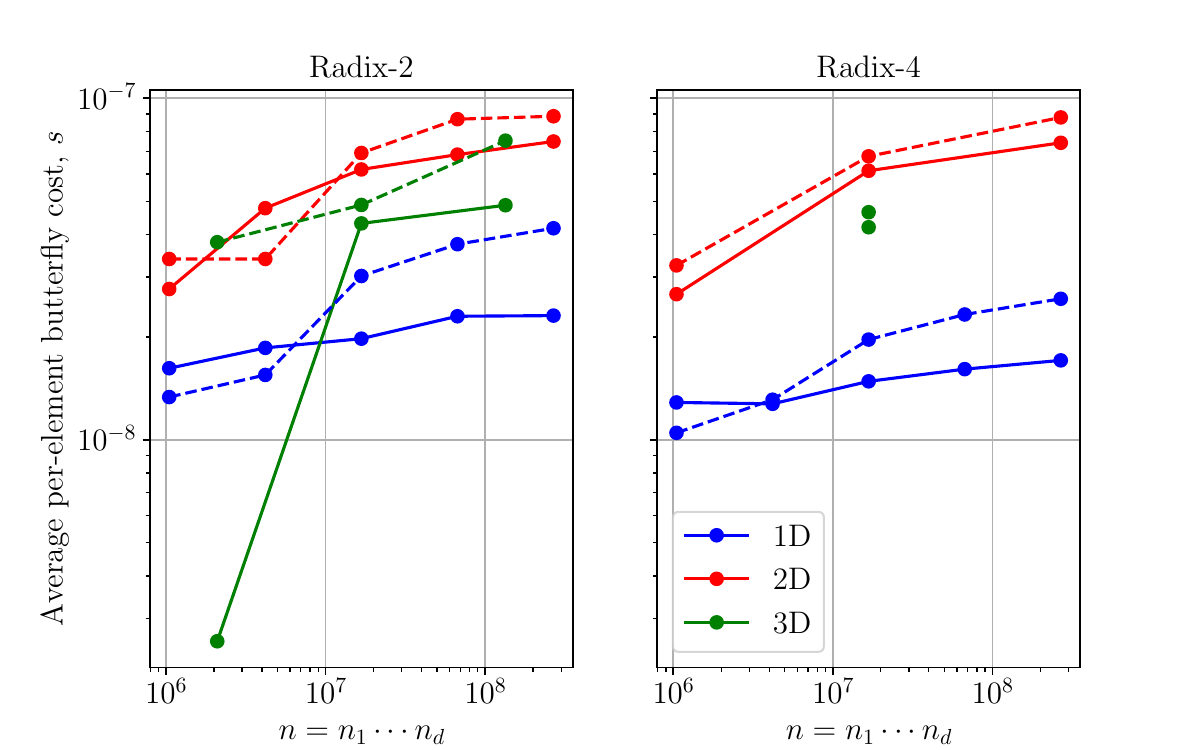}
\caption{Average per-element butterfly cost for computing systems \#1 (solid lines) and \#2 (dashed lines)}
\label{fig:per-element-butterfly-cost}
\end{figure}

Figures~\ref{fig:per-element-permutation-cost} and \ref{fig:per-element-butterfly-cost} provide crucial insights into the scaling behavior of the index-reversal permutation and butterfly operation, respectively:

Average per-element permutation costs increase with total size $n$ but decrease dramatically with dimension $d$ at fixed $n$.
While the average per-element costs reported appear significantly smaller than typical runtime measurement precision, these values are obtained by dividing the total kernel runtime--itself a meaningful and accurate measurement--by the transform size.
The improvement from 1D to 2D and from 2D to 3D each provides nearly an order of magnitude reduction, directly attributable to enhanced spatial locality arising from the Kronecker product structure of permutation matrices.

Per-element butterfly costs show the inverse trend, increasing with both total size $n$ and dimension $d$. 
This reflects the increasing complexity of multi-dimensional memory access patterns required for separable transforms.

These opposing trends suggest that permutation-avoiding techniques provide the greatest benefits for: (i) large 1D problems where permutation costs dominate; (ii) multi-dimensional problems with large sizes in one dimension; and (iii) systems with sufficiently small cache.

The experimental results demonstrate that while theoretical FLOP counts provide useful guidance, real-world FFT performance requires careful consideration of memory hierarchy effects, and the arithmetic intensity of different algorithmic components. 
The permutation-avoiding approach successfully exploits this understanding to achieve meaningful performance improvements in practically relevant scenarios.

\section{Discussion}\label{sec:discussion}
In this section, we investigate the applicability of the permutation-avoiding convolution procedure to other FFT frameworks, namely the mixed-radix case, the prime factor framework, Rader's and Bluestein's algorithms as well as the four-step method.

\subsection{Mixed-radix framework}
Although we aimed to keep our treatment of the Cooley-Tukey framework as general as possible (see Section~\ref{sec:cooley-tukey}), we still assumed that each recursion step subdivides the problem into $r$ equal-sized subproblems. 
That is, we used a uniform radix.
In practice, a more general and widely used approach is the mixed-radix framework, in which $n$ is factored into a product of $t$ integers, i.e., $n = r_1 \cdots r_t$.

The resulting Fourier matrix retains its fundamental properties, as these are independent of the framework adopted. 
However, the factorization itself differs. 
Within the Cooley-Tukey framework, the mixed-radix approach leads to a decomposition of the form
\begin{align}
F_n = A_{\rho,n} P_{\rho,n}^T,
\end{align}
where $\rho := [r_1, \dots, r_t]$. 
Here, $A_{\rho,n} \in \mathbb{C}^{n \times n}$ is still a product of block butterfly matrices, while $P_{\rho,n}$ is a permutation matrix, generally not symmetric, unlike the case in the uniform radix framework discussed in Section~\ref{sec:cooley-tukey}.

The matrix $A_{\rho,n}$ has the form:
\begin{align}
A_{\rho,n} =
B_{r_t,n}
(I_{n/(r_1 \cdots r_{t-1})} \otimes B_{r_{t-1},r_1 \cdots r_{t-1}})
\cdots
(I_{n/(r_1 r_2)} \otimes B_{r_2,r_1 r_2})
(I_{n/r_1} \otimes B_{r_1,r_1}),
\end{align}
where $B_{r,n}$ is defined as in Eq.~\eqref{eq:25}.

Similarly, the permutation matrix is given by:
\begin{align}
P_{\rho,n} =
\Pi_{r_t,n}
(I_{n/(r_1 \cdots r_{t-1})} \otimes \Pi_{r_{t-1},r_1 \cdots r_{t-1}})
\cdots
(I_{n/(r_1 r_2)} \otimes \Pi_{r_2,r_1 r_2})
(I_{n/r_1} \otimes \Pi_{r_1,r_1}),
\end{align}
where $\Pi_{r,n}$ denotes a modulo-$r$ perfect shuffle permutation matrix, that is, the inverse of a modulo-$r$ sort permutation matrix. 
Note that both coincide when $n$ is a power of $r$, as assumed in Section~\ref{sec:cooley-tukey}.

The transpose symmetry of the Fourier matrix still holds:
\begin{align}
F_n = F_n^T = (A_{\rho,n} P_{\rho,n}^T)^T = P_{\rho,n} A_{\rho,n}^T.
\end{align}
Moreover, the inverse Fourier matrix satisfies:
\begin{align}
F_n^{-1} = \overline{F_n}/n = \overline{A_{\rho,n} P_{\rho,n}^T} = \overline{A_{\rho,n}} P_{\rho,n}^T.
\end{align}
Putting all components together in the operator $\mathcal{F}_h$, we obtain:
\begin{align}
\mathcal{F}_h(x)
&= F_n^{-1}((F_n h)\circ (F_n x))\\
&= \overline{A_{\rho,n}} P_{\rho,n}^T((P_{\rho,n} A_{\rho,n}^Th)\circ(P_{\rho,n} A_{\rho,n}^Tx))/ n \nonumber \\
&= \overline{A_{\rho,n}} P_{\rho,n}^TP_{\rho,n}((A_{\rho,n}^Th)\circ(A_{\rho,n}^Tx))/ n \nonumber \\
&= \overline{A_{\rho,n}} ((A_{\rho,n}^Th)\circ(A_{\rho,n}^Tx))/ n \nonumber.
\end{align}
Consequently, the permutation-avoiding expression for $\mathcal{F}_h$ given in Eq.~\eqref{eq:35} remains valid in the mixed-radix setting. 
Furthermore, the same applies to multi-dimensional mixed-radix FFTs, and Algo.~\ref{alg:permutation-avoiding-convolution} remains valid.

\subsection{Prime factor framework}
Prime factor FFTs~\cite{good1958interaction,thomas1963using} rely on a number-theoretic decomposition of the Fourier matrix. 
Suppose $n = n_1 \cdots n_t$ with $n_1, \dots, n_t$ pairwise coprime. 
Then, there exist permutation matrices $\Gamma_\rho, \Upsilon_\rho \in \mathbb{R}^{n \times n}$ such that
\begin{align}
F_n = \Gamma_\rho (F_{n_t} \otimes \cdots \otimes F_{n_1}) \Upsilon_\rho^T.
\end{align}
Assume now that each Fourier matrix $F_{n_q}$ can be locally factored as $
F_{n_q}=A_{n_q}P_{n_q}^T$ for $q=1,\dots,t$, with $A_{n_q}\in\mathbb{C}^{n_q \times n_q}$ and $P_{n_q}\in\mathbb{R}^{n_q \times n_q}$. 
Then, we have
\begin{align}
F_n
=&\,F_n^T\\
=&\,\left(\Gamma_\rho(F_{n_t}\otimes\cdots\otimes F_{n_1})\Upsilon_\rho^T\right)^T\nonumber\\
=&\,\left(\Gamma_\rho(A_{n_t}P_{n_t}^T)\otimes\cdots\otimes(A_{n_1}P_{n_1}^T)\Upsilon_\rho^T\right)^T\nonumber\\
=&\,\Upsilon_\rho (P_{n_t} \otimes \cdots \otimes P_{n_1})(A_{n_t}^T \otimes \cdots \otimes A_{n_1}^T) \Gamma_\rho^T\nonumber
\end{align}
and, similarly for the inverse Fourier matrix:
\begin{align}
F_n^{-1}
=&\,\overline{F_n}/n\\
=&\,\overline{\Gamma_\rho(F_{n_t}\otimes\cdots\otimes F_{n_1})\Upsilon_\rho^T}/n\nonumber\\
=&\,\overline{\Gamma_\rho(A_{n_t}P_{n_t}^T)\otimes\cdots\otimes(A_{n_1}P_{n_1}^T)\Upsilon_\rho^T}/n\nonumber\\
=&\,\Gamma_\rho(\overline{A_{n_t}}\otimes\cdots\otimes\overline{A_{n_1}})(P_{n_t}^T \otimes\cdots\otimes P_{n_1}^T)\Upsilon_\rho^T/n\nonumber
\end{align}
so that the convolution $\mathcal{F}_h(x) = F_n^{-1}((F_nh)\circ(F_n x))$ becomes
\begin{align}
\mathcal{F}_h(x)=
\Gamma_\rho (\overline{A_{n_t}}\otimes\cdots\otimes\overline{A_{n_1}})
\left[
\left((A_{n_t}^T\otimes\cdots\otimes A_{n_1}^T)\Gamma_\rho^Th\right)
\circ
\left((A_{n_t}^T\otimes\cdots\otimes A_{n_1}^T)\Gamma_\rho^Tx\right)
\right]/n.
\end{align}
This formulation demonstrates that, the prime factor FFT-based convolution can be performed using only three permutations, $\Gamma_\rho^T$ for the forward transforms, and $\Gamma_\rho$ for the backward transform, instead of six, thus enabling some acceleration through permutation-avoidance.
Note that these observations generalize directly to multi-dimensional prime factor FFTs. 

\subsection{Rader's and Bluestein's frameworks}
Both Rader’s FFT~\cite{rader1968discrete} and the Bluestein algorithm~\cite{bluestein1970linear} reformulate the discrete Fourier transform (DFT) as a convolution, which can then be efficiently computed using any suitable FFT framework. 
If the chosen FFT framework supports permutation-avoiding techniques, as demonstrated earlier, then Rader’s and Bluestein’s methods can also benefit from these optimizations.

\subsection{Four-step framework}
One of the most performant FFT frameworks is a blocking strategy initially proposed in~\cite{gentleman1966fast}. 
After being revisited in~\cite{bailey1990ffts}, it became known as the four-step FFT, and is sometimes even referred to as Bailey's algorithm.

The starting point of the four-step FFT is the radix-$r$ splitting. 
That is, suppose $n = rm$, then we have
\begin{align}
F_n
=
(F_r \otimes I_m)\,\text{diag}(I_m, \Omega_{r,m}, \dots, \Omega_{r,m}^{r-1})(I_r \otimes F_m)\Pi_{r,n}^T.
\end{align}
Let $x \in \mathbb{C}^n$, and denote by $x_{m \times r}$ its $m \times r$ column-major matrixification. Then, for any $A \in \mathbb{C}^{r \times r}$, we have $((A \otimes I_m)x)_{m \times r} = x_{m \times r}A^T$.
Applying this identity yields
\begin{align}
(F_n x)_{m \times r}
&= ((F_r \otimes I_m)\,\text{diag}(I_m, \Omega_{r,m}, \dots, \Omega_{r,m}^{r-1})(I_r \otimes F_m)\Pi_{r,n}^T x)_{m \times r} \\
&= (\text{diag}(I_m, \Omega_{r,m}, \dots, \Omega_{r,m}^{r-1})(I_r \otimes F_m)\Pi_{r,n}^T x)_{m \times r} F_r \nonumber \\
&= \left((F_n)_{1:m,1:r} \circ ((I_r \otimes F_m)\Pi_{r,n}^T x)_{m \times r}\right) F_r. \nonumber
\end{align}
Using the additional Kronecker identity $((I_r \otimes A)x)_{m \times r} = A x_{m \times r}$ for $A \in \mathbb{C}^{m \times m}$, we obtain
\begin{align}
(F_n x)_{m \times r} =
\left((F_n)_{1:m,1:r} \circ (F_m (\Pi_{r,n}^T x)_{m \times r})\right) F_r.
\end{align}
Since $\Pi_{r,n}^T$ is a modulo-$r$ sort permutation matrix, we have $(\Pi_{r,n}^T x)_{m \times r} = x_{r \times m}^T$, leading to
\begin{align}\label{eq:100}
(F_n x)_{m \times r} =
\left((F_n)_{1:m,1:r} \circ (F_m x_{r \times m}^T)\right) F_r.
\end{align}
This is the final form of the four-step FFT. 
Note that the single-vector DFT $F_n x$ is transformed into a sequence of DFTs on multiple vectors.

If the constituent DFTs are implemented via a radix-2 Cooley-Tukey algorithm, the total number of FLOP required by Eq.~\eqref{eq:100} is $5n \log_2 n$, which matches the complexity of the radix-2 Cooley-Tukey FFT. 
However, while the four-step framework does not reduce the overall FLOP count, its blocking strategy dramatically improves data locality and reuse, leading to some of the fastest known implementations.

Applying the four-step approach to the inverse DFT yields
\begin{align}
(F_n^{-1} y)_{r \times m}
&= (\overline{F_n} y)_{r \times m} / n \\
&= \left(\overline{(F_n)_{1:r,1:m}} \circ (\overline{F_r} y_{m \times r}^T)\right) \overline{F_m} / n \\
&= n\left(\overline{(F_n)_{1:r,1:m}} \circ (F_r^{-1} y_{m \times r}^T)\right) F_m^{-1}.
\end{align}
Therefore, the convolution $\mathcal{F}_h(x) = F_n^{-1}((F_nh) \circ (F_n x))$ can be written as
\begin{align}
(\mathcal{F}_h(x))_{r \times m}
&= n \left(\overline{(F_n^{-1})_{1:r,1:m}} \circ \left(F_r^{-1} 
\left[(F_n h)_{m \times r} \circ (F_n x)_{m \times r}\right]^T
\right)\right) F_m^{-1}.
\end{align}
Upon detailed examination of the four-step framework and its application to discrete convolution, we conclude that it is not amenable to permutation-avoiding techniques.

\section{Conclusion}\label{sec:conclusion}
This work presents a practical analysis of permutation-avoiding convolution algorithms within the general radix Cooley-Tukey framework, extending from one-dimensional to multi-dimensional transforms. 
In doing so, we demonstrate achievable performance improvements for the computation of discrete convolutions.

Our investigation reveals that index-reversal permutations, while representing zero arithmetic intensity, can dominate FFT execution time due to their memory-bound nature and poor cache behavior. 
By reformulating the mathematical structure of FFT-based convolution, as was previously done in~\cite{van1992computational}, we show that these permutations can be eliminated in Cooley-Tukey frameworks, as well as in other types of implementations, in 1D and multiple dimensions. 
In the case of implementations based on the Cooley-Tukey framework, this transformation eliminates the overhead of three permutation applications per convolution, limiting the cost of discrete convolution to butterfly operations and component-wise multiplication with the transformed filter.

The experimental validation across multiple radices (2, 4, and 8) and dimensions (1D, 2D, and 3D) provides several insights. 
First, theoretical FLOP counts alone provide insufficient guidance for real-world performance optimization---our radix-4 implementation consistently outperformed radix-8 despite the latter's 20\% FLOP advantage, highlighting the critical importance of memory access patterns and register utilization. 
Second, the effectiveness of permutation-avoiding techniques exhibits strong dimensional dependence: 1D problems show dramatic average speedups (up to 4.04X for sufficiently large transforms), while multi-dimensional cases demonstrate more modest but consistent improvements due to improved spatial locality in higher-dimensional Kronecker products of products, and increased cost of multi-dimensional butterfly operations.

Our scaling analysis reveals the fundamental trade-off between permutation and butterfly costs across dimensions. 
While per-element permutation costs decrease dramatically with dimension due to enhanced spatial locality, butterfly operations become increasingly expensive, due to the need to fetch noncontiguous data from one dimension to another or, alternatively, apply transpositions.
This analysis provides clear guidance for when permutation-avoiding approaches offer maximum benefit: large 1D problems where permutation costs dominate, multi-dimensional problems with large sizes in one dimension, and systems with sufficiently small cache.

The practical impact extends beyond our specific implementation. 
We demonstrate that the permutation-avoiding approach generalizes to mixed-radix, prime factor, and other Cooley-Tukey variants, making it applicable to production FFT libraries. 
Our modular implementation architecture, which cleanly separates permutation and butterfly operations, enables precise performance analysis and could inform the design of next-generation discrete convolution kernels.

Several directions for future work emerge from this analysis. 
First, the integration of permutation-avoiding techniques with advanced optimization strategies such as cache blocking, vectorization, and multi-threading remains unexplored. 
In particular, our multi-dimensional implementations could benefit significantly from blocking strategies that process multiple rows or columns simultaneously, improving cache locality and reducing the gather-scatter overhead that currently limits 2D and 3D performance.
Multi-dimensional data transpositions can also be made more efficient.
Second, the extension to mixed-radix and other FFT frameworks could broaden applicability. 
Third, the potential for hardware-specific optimizations, particularly for architectures with different cache hierarchies or specialized memory access patterns, warrants investigation.
Finally, it could be interesting to compare permutation-avoiding implementations with the type of code optimization achievable using the SPIRAL framework~\cite{franchetti_spiral_2018}.

Our results suggest that developers of high-performance FFT libraries should consider supporting permutation-avoiding convolution kernels. 
The performance improvements demonstrated across different problem sizes and dimensions, combined with the applicability to existing FFT frameworks, make this a compelling addition to the computational toolkit for signal processing, scientific computing, and machine learning applications.

This work reinforces the importance of considering memory hierarchy effects alongside computational complexity in algorithm design. 
As the gap between processor and memory performance continues to widen, techniques that minimize memory movement while preserving computational correctness become increasingly valuable for achieving optimal performance in numerical computing.

\section{Acknowledgment}
We would like to thank Erik Wernersson from the Karolinska Instutet, who brought to our attention the fact that the permutation matrices applied to the convolution filter can be factored out along with the permutation of the forward data transformation, thus leading to the complete cancellation of three permutations from the convolution operation, i.e., see Eq.~\eqref{eq:34}, which is referred to as \textit{fast convolution} in \cite{van1992computational}.
Before Erik Wernersson's comment, an earlier version of this document did not factor out the permutation from the forward filter transformation, which led to a strategy recommended for use in cases where multiple convolutions are applied with a single filter, permuted once and for all, in an offline stage, whereas the online stage of the strategy remained permutation-free.

\bibliographystyle{ACM-Reference-Format}
\bibliography{references}

\appendix
\section{FFT and convolution sub-kernel algorithms}\label{app}
In the algorithms presented in Section~\ref{sec:algorithms}, for any given $x\in\mathbb{C}^n$, we need to evaluate $A_{r,n}x$, $\overline{A_{r,n}}x$ and $A_{r,n}^Tx$.
Irrespective of the radix, these kernels are implemented as described in Algos.~\ref{alg:butterfly-related-kernel} and \ref{alg:transposed-butterfly-kernel}.
In what follows, we show detailed implementations of these kernels for radices 2, 4, and 8.

\begin{algorithm}[H]
\caption{Butterfly-related kernels}
\label{alg:butterfly-related-kernel}
\begin{algorithmic}[1]
\REQUIRE{$x\in\mathbb{C}^n$}
\ENSURE{$x:=A_{r,n}x$ or $x:=\overline{A_{r,n}}x$}
\FOR{$q=1,\dots,t$}
\STATE{$k:=r^q$}
\STATE{$n_b:=n/k$}
%\STATE{$\ell:=k/r$}
\FOR{$b=1,\dots,n_b$}
\STATE{\hspace{.4cm}$x_{(b-1)k+1:bk}:=B_{r,k}x_{(b-1)k+1:bk}$}
\STATE{or $x_{(b-1)k+1:bk}:=\overline{B_{r,k}}x_{(b-1)k+1:bk}$}
\ENDFOR
\ENDFOR
\RETURN{$x$}
\end{algorithmic}
\end{algorithm}

\begin{algorithm}[H]
\caption{Transposed butterfly kernels}
\label{alg:transposed-butterfly-kernel}
\begin{algorithmic}[1]
\REQUIRE{$x\in\mathbb{C}^n$}
\ENSURE{$x:=A_{r,n}^Tx$}
\FOR{$q=t,\dots,1$}
\STATE{$k:=r^q$}
\STATE{$n_b:=n/k$}
%\STATE{$\ell:=k/r$}
\FOR{$b=1,\dots,n_b$}
\STATE{$x_{(b-1)k+1:bk}:=B_{r,k}^Tx_{(b-1)k+1:bk}$}
\ENDFOR
\ENDFOR
\RETURN{$x$}
\end{algorithmic}
\end{algorithm}

\subsection{Radix-2 butterfly-related kernels}
For the radix-2 case, we have
\begin{align}
B_{2,k}=
(F_2\otimes I_{k/2})\,\text{diag}\left(I_{k/2},\Omega_{2,k/2}\right)
\end{align}
where 
$F_2=\begin{bmatrix}1& 1\\1&-1\end{bmatrix}$ 
so that
$B_{2,k}=
\begin{bmatrix}
I_{k/2} &  \Omega_{2,k/2}\\
I_{k/2} & -\Omega_{2,k/2}
\end{bmatrix}$.\\
Then, for all $x\in\mathbb{C}^k$, as we let $z_i:=x_{(i-1)k/2+1:ik/2}$ for $i=1,2$, we obtain
\begin{align}
B_{2,k}x=&\,
\begin{bmatrix}
z_{1}+\Omega_{2,k/2}z_{2}\\
z_{1}-\Omega_{2,k/2}z_{2}\\
\end{bmatrix}.
\end{align}
We then let $\tau:=\Omega_{2,k/2}z_{2}$ so that
\begin{align}
B_{2,k}x=&\,
\begin{bmatrix}
z_{1}+\tau\\
z_{1}-\tau
\end{bmatrix}.
\end{align}
This leads to Algo.~\ref{alg:butterfly-kernel-radix-2} for the computation of $x\mapsto A_{2,n}x$.
\begin{algorithm}[H]
\caption{Butterfly kernel of radix-2}
\label{alg:butterfly-kernel-radix-2}
\begin{algorithmic}[1]
\REQUIRE{$x\in\mathbb{C}^n$, $n=2^t$}
\ENSURE{$x:=A_{2,n}x$}
\FOR{$q=1,\dots,t$}
\STATE{$k:=2^q$}
\STATE{$n_b:=n/k$}
\STATE{$\ell:=k/2$}
\FOR{$b=1,\dots,n_b$}
\FOR{$j=1\dots,\ell$}
\STATE{$\tau:=\omega_k^{j-1}x_{(b-1)k+j+\ell}$}
\STATE{$x_{(b-1)k+j+\ell}:=x_{(b-1)k+j}-\tau$}
\STATE{$x_{(b-1)k+j}:=x_{(b-1)k+j}+\tau$}
\ENDFOR
\ENDFOR
\ENDFOR
\RETURN{$x$}
\end{algorithmic}
\end{algorithm}
The radix-2 conjugate butterfly kernel relies on the expression
\begin{align}
\overline{B_{2,k}}x=&\,
\begin{bmatrix}
z_{1}+\overline{\Omega_{2,k/2}}z_{2}\\
z_{1}-\overline{\Omega_{2,k/2}}z_{2}\\
\end{bmatrix}
\end{align}
so that, as we let $\tau:=\overline{\Omega_{2,k/2}}z_{2}$, we still have
\begin{align}
\overline{B_{2,k}}x=&\,
\begin{bmatrix}
z_{1}+\tau\\
z_{1}-\tau
\end{bmatrix}.
\end{align}
This leads to Algo.~\ref{alg:conjugate-butterfly-kernel-radix-2} for the computation of $x\mapsto\overline{A_{2,n}}x$.
\begin{algorithm}[H]
\caption{Conjugate butterfly kernel of radix-2}
\label{alg:conjugate-butterfly-kernel-radix-2}
\begin{algorithmic}[1]
\REQUIRE{$x\in\mathbb{C}^n$, $n=2^t$}
\ENSURE{$x:=\overline{A_{2,n}}x$}
\FOR{$q=1,\dots,t$}
\STATE{$k:=2^q$}
\STATE{$n_b:=n/k$}
\STATE{$\ell:=k/2$}
\FOR{$b=1,\dots,n_b$}
\FOR{$j=1\dots,\ell$}
\STATE{$\tau:=\overline{\omega_k^{j-1}}x_{(b-1)k+j+\ell}$}
\STATE{$x_{(b-1)k+j+\ell}:=x_{(b-1)k+j}-\tau$}
\STATE{$x_{(b-1)k+j}:=x_{(b-1)k+j}+\tau$}
\ENDFOR
\ENDFOR
\ENDFOR
\RETURN{$x$}
\end{algorithmic}
\end{algorithm}
On the other hand, we also have
\begin{align}
B_{2,k}^Tx=
\begin{bmatrix}
z_{1}+z_{2}\hfill\\
\Omega_{2,k/2}z_{1}-\Omega_{2,k/2}z_{2}
\end{bmatrix}
\end{align}
so that if we let $\tau:=z_2$, we obtain
\begin{align}
B_{2,k}^Tx=
\begin{bmatrix}
z_{1}+\tau\hfill\\
\Omega_{2,k/2}(z_{1}-\tau)
\end{bmatrix}.
\end{align}
This leads to Algo.~\ref{alg:transposed-butterfly-kernel-radix-2} for the computation of $x\mapsto A_{2,n}^Tx$.
\begin{algorithm}[H]
\caption{Transposed butterfly kernel of radix-2}
\label{alg:transposed-butterfly-kernel-radix-2}
\begin{algorithmic}[1]
\REQUIRE{$x\in\mathbb{C}^n$, $n=2^t$}
\ENSURE{$x:=A_{2,n}^Tx$}
\FOR{$q=t,\dots,1$}
\STATE{$k:=2^q$}
\STATE{$n_b:=n/k$}
\STATE{$\ell:=k/2$}
\FOR{$b=1,\dots,n_b$}
\FOR{$j=1\dots,\ell$}
\STATE{$\tau:=x_{(b-1)k+j+\ell}$}
\STATE{$x_{(b-1)k+j+\ell}:=\omega_k^{j-1}(x_{(b-1)k+j}-\tau)$}
\STATE{$x_{(b-1)k+j}:=x_{(b-1)k+j}+\tau$}
\ENDFOR
\ENDFOR
\ENDFOR
\RETURN{$x$}
\end{algorithmic}
\end{algorithm}

\subsection{Radix-4 butterfly-related kernels}
For the radix-4 case, we have
\begin{align}
B_{4,k}=
(F_4\otimes I_{k/4})\,\text{diag}\left(I_{k/4},\Omega_{4,k/4},\Omega_{4,k/4}^2,\Omega_{4,k/4}^3\right)
\end{align}
where
\begin{align}
F_4=
\begin{bmatrix}
1& 1& 1& 1\\
1&-i&-1& i\\
1&-1& 1&-1\\
1& i&-1&-i
\end{bmatrix}
\end{align}
so that
\begin{align}
B_{4,k}=
\begin{bmatrix}
I_{k/4}&       \Omega_{4,k/4}& \Omega_{4,k/4}^2&       \Omega_{4,k/4}^3\\
I_{k/4}&-i\cdot\Omega_{4,k/4}&-\Omega_{4,k/4}^2& i\cdot\Omega_{4,k/4}^3\\
I_{k/4}&      -\Omega_{4,k/4}& \Omega_{4,k/4}^2&      -\Omega_{4,k/4}^3\\
I_{k/4}& i\cdot\Omega_{4,k/4}&-\Omega_{4,k/4}^2&-i\cdot\Omega_{4,k/4}^3
\end{bmatrix}.
\end{align}
Then, for all $x\in\mathbb{C}^k$, as we let $z_i:=x_{(i-1)k/4+1:ik/4}$ for $i=1,\dots,4$, we obtain
\begin{align}
B_{4,k}x=&\,
\begin{bmatrix}
z_{1}+\Omega_{4,k/4}z_{2}+\Omega_{4,k/4}^2z_{3}+\Omega_{4,k/4}^3z_{4}\hfill\\
z_{1}-i\cdot\Omega_{4,k/4}z_{2}-\Omega_{4,k/4}^2z_{3}+i\cdot\Omega_{4,k/4}^3z_{4}\\
z_{1}-\Omega_{4,k/4}z_{2}+\Omega_{4,k/4}^2z_{3}-\Omega_{4,k/4}^3z_{4}\hfill\\
z_{1}+i\cdot\Omega_{4,k/4}z_{2}-\Omega_{4,k/4}^2z_{3}-i\cdot\Omega_{4,k/4}^3z_{4}\hfill
\end{bmatrix}.
\end{align}
Once we introduce
\begin{align}
&\tau_1:=z_1+\Omega_{4,k/4}^2z_3,\hspace{1.02cm}
\tau_2:=z_1-\Omega_{4,k/4}^2z_3\\
&\tau_3:=\Omega_{4,k/4}z_2+\Omega_{4,k/4}^3z_4,\;\;
\tau_4:=\Omega_{4,k/4}z_2-\Omega_{4,k/4}^3z_4\nonumber
\end{align}
we get
\begin{align}
B_{4,k}x=&\,
\begin{bmatrix}
\tau_1+\tau_3\hfill\\
\tau_2-i\cdot\tau_4\hfill\\
\tau_1-\tau_3\hfill\\
\tau_2+i\cdot\tau_4\hfill
\end{bmatrix}.
\end{align}
This leads to Algo.~\ref{alg:butterfly-kernel-radix-4} for the computation of $x\mapsto A_{4,n}x$.
\begin{algorithm}[H]
\caption{Butterfly kernel of radix-4}
\label{alg:butterfly-kernel-radix-4}
\begin{algorithmic}[1]
\REQUIRE{$x\in\mathbb{C}^n$, $n=4^t$}
\ENSURE{$x:=A_{4,n}x$}
\FOR{$q=1,\dots,t$}
\STATE{$k:=4^q$}
\STATE{$n_b:=n/k$}
\STATE{$\ell:=k/4$}
\FOR{$b=1,\dots,n_b$}
\FOR{$j=1\dots,\ell$}
\STATE{$z_1:=x_{(b-1)k+j}$}
\STATE{$z_2:=\omega_k^{j-1}x_{(b-1)k+\ell+j}$}
\STATE{$z_3:=\omega_k^{2(j-1)}x_{(b-1)k+2\ell+j}$}
\STATE{$z_4:=\omega_k^{3(j-1)}x_{(b-1)k+3\ell+j}$}
\STATE{$\tau_1:=z_1+z_3$, $\tau_2:=z_1-z_3$}
\STATE{$\tau_3:=z_2+z_4$, $\tau_4:=z_2-z_4$}
\STATE{$x_{(b-1)k+j}:=\tau_1+\tau_3$}
\STATE{$x_{(b-1)k+\ell+j}:=\tau_2-i\cdot \tau_4$}
\STATE{$x_{(b-1)k+2\ell+j}:=\tau_1-\tau_3$}
\STATE{$x_{(b-1)k+3\ell+j}:=\tau_2+i\cdot \tau_4$}
\ENDFOR
\ENDFOR
\ENDFOR
\RETURN{$x$}
\end{algorithmic}
\end{algorithm}
The radix-4 conjugate butterfly kernel relies on the expression
\begin{align}
\overline{B_{4,k}}x=&\,
\begin{bmatrix}
z_{1}+\overline{\Omega_{4,k/4}}z_{2}+\overline{\Omega_{4,k/4}^2}z_{3}+\overline{\Omega_{4,k/4}^3}z_{4}\hfill\\
z_{1}+i\cdot\overline{\Omega_{4,k/4}}z_{2}-\overline{\Omega_{4,k/4}^2}z_{3}-i\cdot\overline{\Omega_{4,k/4}^3}z_{4}\\
z_{1}-\overline{\Omega_{4,k/4}}z_{2}+\overline{\Omega_{4,k/4}^2}z_{3}-\overline{\Omega_{4,k/4}^3}z_{4}\hfill\\
z_{1}-i\cdot\overline{\Omega_{4,k/4}}z_{2}-\overline{\Omega_{4,k/4}^2}z_{3}+i\cdot\overline{\Omega_{4,k/4}^3}z_{4}\hfill
\end{bmatrix}
\end{align}
so that, as we introduce
\begin{align}
&\tau_1:=z_1+\overline{\Omega_{4,k/4}^2}z_3,\hspace{1.02cm}
\tau_2:=z_1-\overline{\Omega_{4,k/4}^2}z_3\\
&\tau_3:=\overline{\Omega_{4,k/4}}z_2+\overline{\Omega_{4,k/4}^3}z_4,\;\;
\tau_4:=\overline{\Omega_{4,k/4}}z_2-\overline{\Omega_{4,k/4}^3}z_4\nonumber
\end{align}
we get
\begin{align}
\overline{B_{4,k}}x=&\,
\begin{bmatrix}
\tau_1+\tau_3\hfill\\
\tau_2+i\cdot\tau_4\hfill\\
\tau_1-\tau_3\hfill\\
\tau_2-i\cdot\tau_4\hfill
\end{bmatrix}.
\end{align}
This leads to Algo.~\ref{alg:conjugate-butterfly-kernel-radix-4} for the computation of $x\mapsto \overline{A_{4,n}}x$.
\begin{algorithm}[H]
\caption{Conjugate butterfly kernel of radix-4}
\label{alg:conjugate-butterfly-kernel-radix-4}
\begin{algorithmic}[1]
\REQUIRE{$x\in\mathbb{C}^n$, $n=4^t$}
\ENSURE{$x:=\overline{A_{4,n}}x$}
\FOR{$q=1,\dots,t$}
\STATE{$k:=4^q$}
\STATE{$n_b:=n/k$}
\STATE{$\ell:=k/4$}
\FOR{$b=1,\dots,n_b$}
\FOR{$j=1\dots,\ell$}
\STATE{$z_1:=x_{(b-1)k+j}$}
\STATE{$z_2:=\overline{\omega_k^{j-1}}x_{(b-1)k+\ell+j}$}
\STATE{$z_3:=\overline{\omega_k^{2(j-1)}}x_{(b-1)k+2\ell+j}$}
\STATE{$z_4:=\overline{\omega_k^{3(j-1)}}x_{(b-1)k+3\ell+j}$}
\STATE{$\tau_1:=z_1+z_3$, $\tau_2:=z_1-z_3$}
\STATE{$\tau_3:=z_2+z_4$, $\tau_4:=z_2-z_4$}
\STATE{$x_{(b-1)k+j}:=\tau_1+\tau_3$}
\STATE{$x_{(b-1)k+\ell+j}:=\tau_2+i\cdot \tau_4$}
\STATE{$x_{(b-1)k+2\ell+j}:=\tau_1-\tau_3$}
\STATE{$x_{(b-1)k+3\ell+j}:=\tau_2-i\cdot \tau_4$}
\ENDFOR
\ENDFOR
\ENDFOR
\RETURN{$x$}
\end{algorithmic}
\end{algorithm}
The radix-4 transposed butterfly kernel relies on the expression
\begin{align}
B_{4,k}^Tx=&\,
\begin{bmatrix}
z_{1}+z_{2}+z_{3}+z_{4}\hfill\\
\Omega_{4,k/4}z_{1}-i\cdot\Omega_{4,k/4}z_{2}-\Omega_{4,k/4}z_{3}+i\cdot\Omega_{4,k/4}z_{4}\\
\Omega_{4,k/4}^2z_{1}-\Omega_{4,k/4}^2z_{2}+\Omega_{4,k/4}^2z_{3}-\Omega_{4,k/4}^2z_{4}\hfill\\
\Omega_{4,k/4}^3z_{1}+i\cdot\Omega_{4,k/4}^3z_{2}-\Omega_{4,k/4}^3z_{3}-i\cdot\Omega_{4,k/4}^3z_{4}\hfill
\end{bmatrix}.
\end{align}
As we introduce
\begin{align}
\tau_1:=z_1+z_3,\;\;
\tau_2:=\Omega_{4,k/4}(z_1-z_3)\\
\tau_3:=z_2+z_4,\;\;
\tau_4:=\Omega_{4,k/4}(z_2-z_4)\nonumber
\end{align}
we obtain
\begin{align}
B_{4,k}^Tx=&\,
\begin{bmatrix}
\tau_1+\tau_3\hfill\\
\tau_2-i\cdot\tau_4\hfill\\
\Omega_{4,k/4}^2(\tau_1-\tau_3)\hfill\\
\Omega_{4,k/4}^2(\tau_2+i\cdot\tau_4)\hfill
\end{bmatrix}.
\end{align}
This leads to Algo.~\ref{alg:transposed-butterfly-kernel-radix-4} for the computation of $x\mapsto A_{4,n}^Tx$.
\begin{algorithm}[H]
\caption{Transposed butterfly kernel of radix-4}
\label{alg:transposed-butterfly-kernel-radix-4}
\begin{algorithmic}[1]
\REQUIRE{$x\in\mathbb{C}^n$, $n=4^t$}
\ENSURE{$x:=A_{4,n}^Tx$}
\FOR{$q=t,\dots,1$}
\STATE{$k:=4^q$}
\STATE{$n_b:=n/k$}
\STATE{$\ell:=k/4$}
\FOR{$b=1,\dots,n_b$}
\FOR{$j=1\dots,\ell$}
\STATE{$z_1:=x_{(b-1)k+j}$}
\STATE{$z_2:=x_{(b-1)k+\ell+j}$}
\STATE{$z_3:=x_{(b-1)k+2\ell+j}$}
\STATE{$z_4:=x_{(b-1)k+3\ell+j}$}
\STATE{$\tau_1:=z_1+z_3$, $\tau_2:=\omega_{k}^{j-1}(z_1-z_3)$}
\STATE{$\tau_3:=z_2+z_4$, $\tau_4:=\omega_{k}^{j-1}(z_2-z_4)$}
\STATE{$x_{(b-1)k+j}:=\tau_1+\tau_3$}
\STATE{$x_{(b-1)k+\ell+j}:=\tau_2-i\cdot \tau_4$}
\STATE{$x_{(b-1)k+2\ell+j}:=\omega_k^{2(j-1)}(\tau_1-\tau_3)$}
\STATE{$x_{(b-1)k+3\ell+j}:=\omega_k^{2(j-1)}(\tau_2+i\cdot\tau_4)$}
\ENDFOR
\ENDFOR
\ENDFOR
\RETURN{$x$}
\end{algorithmic}
\end{algorithm}

\subsection{Radix-8 butterfly-related kernels}
For the radix-8 case, we have
\begin{align}
B_{8,k}=
(F_8\otimes I_{k/8})\,\text{diag}\left(I_{k/8},\Omega_{8,k/8},\Omega_{8,k/8}^2,\Omega_{8,k/8}^3,\Omega_{8,k/8}^4,\Omega_{8,k/8}^5,\Omega_{8,k/8}^6,\Omega_{8,k/8}^7\right)
\end{align}
where
\begin{align}
F_8=
\begin{bmatrix}
1& 1& 1& 1& 1& 1& 1& 1\\
1& a&-i& b&-1&-a& i&-b\\
1&-i&-1& i& 1&-i&-1& i\\
1& b& i& a&-1&-b&-i&-a\\
1&-1& 1&-1& 1&-1& 1&-1\\
1&-a&-i&-b&-1& a& i& b\\
1& i&-1&-i& 1& i&-1&-i\\
1&-b& i&-a&-1& b&-i& a
\end{bmatrix}
\end{align}
in which $a=(1-i)/\sqrt{2}$ and $b=-(1+i)/\sqrt{2}=-ia$, so that
\begin{align}
\hspace{-.5cm}
B_{8,k}=
\begin{bmatrix}
I_{k/8}&       \Omega_{8,k/8}&       \Omega_{8,k/8}^2& 
      \Omega_{8,k/8}^3& \Omega_{8,k/8}^4&        \Omega_{8,k/8}^5&
       \Omega_{8,k/8}^6&       \Omega_{8,k/8}^7\\
I_{k/8}& a\cdot\Omega_{8,k/8}&-i\cdot\Omega_{8,k/8}^2& 
 b\cdot\Omega_{8,k/8}^3&-\Omega_{8,k/8}^4&-a\cdot\Omega_{8,k/8}^5&
 i\cdot\Omega_{8,k/8}^6&-b\cdot\Omega_{8,k/8}^7\\
I_{k/8}&-i\cdot\Omega_{8,k/8}&      -\Omega_{8,k/8}^2& 
 i\cdot\Omega_{8,k/8}^3& \Omega_{8,k/8}^4&-i\cdot\Omega_{8,k/8}^5&
      -\Omega_{8,k/8}^6& i\cdot\Omega_{8,k/8}^7\\
I_{k/8}& b\cdot\Omega_{8,k/8}& i\cdot\Omega_{8,k/8}^2& 
 a\cdot\Omega_{8,k/8}^3&-\Omega_{8,k/8}^4&-b\cdot\Omega_{8,k/8}^5&
-i\cdot\Omega_{8,k/8}^6&-a\cdot\Omega_{8,k/8}^7\\
I_{k/8}&      -\Omega_{8,k/8}&       \Omega_{8,k/8}^2&
      -\Omega_{8,k/8}^3& \Omega_{8,k/8}^4&      -\Omega_{8,k/8}^5& 
       \Omega_{8,k/8}^6&      -\Omega_{8,k/8}^7\\
I_{k/8}&-a\cdot\Omega_{8,k/8}&-i\cdot\Omega_{8,k/8}^2&
-b\cdot\Omega_{8,k/8}^3&-\Omega_{8,k/8}^4& a\cdot\Omega_{8,k/8}^5& 
 i\cdot\Omega_{8,k/8}^6& b\cdot\Omega_{8,k/8}^7\\
I_{k/8}& i\cdot\Omega_{8,k/8}&      -\Omega_{8,k/8}^2&
-i\cdot\Omega_{8,k/8}^3& \Omega_{8,k/8}^4& i\cdot\Omega_{8,k/8}^5&
      -\Omega_{8,k/8}^6&-i\cdot\Omega_{8,k/8}^7\\
I_{k/8}&-b\cdot\Omega_{8,k/8}& i\cdot\Omega_{8,k/8}^2&
-a\cdot\Omega_{8,k/8}^3&-\Omega_{8,k/8}^4& b\cdot\Omega_{8,k/8}^5&
-i\cdot\Omega_{8,k/8}^6& a\cdot\Omega_{8,k/8}^7
\end{bmatrix}.
\end{align}
Then, for all $x\in\mathbb{C}^k$, as we let $z_i:=x_{(i-1)k/8+1:ik/8}$ for $i=1,\dots,8$, we have
\begin{align}
\hspace{-.4cm}
B_{8,k}x=&\,
\begin{bmatrix}
z_{1}
+\Omega_{8,k/8}z_{2}
+\Omega_{8,k/8}^2z_{3}
+\Omega_{8,k/8}^3z_{4}
+\Omega_{8,k/8}^4z_{5}
+\Omega_{8,k/8}^5z_{6}
+\Omega_{8,k/8}^6z_{7}
+\Omega_{8,k/8}^7z_{8}\hfill\\
z_{1}
+a\cdot\Omega_{8,k/8}z_{2}
-i\cdot\Omega_{8,k/8}^2z_{3}
+b\cdot\Omega_{8,k/8}^3z_{4}
-\Omega_{8,k/8}^4z_{5}
-a\cdot\Omega_{8,k/8}^5z_{6}
+i\cdot\Omega_{8,k/8}^6z_{7}
-b\cdot\Omega_{8,k/8}^7z_{8}\hfill\\
z_{1}
-i\cdot\Omega_{8,k/8}z_{2}
-\Omega_{8,k/8}^2z_{3}
+i\cdot\Omega_{8,k/8}^3z_{4}
+\Omega_{8,k/8}^4z_{5}
-i\cdot\Omega_{8,k/8}^5z_{6}
-\Omega_{8,k/8}^6z_{7}
+i\cdot\Omega_{8,k/8}^7z_{8}\hfill\\
z_{1}
+b\cdot\Omega_{8,k/8}z_{2}
+i\cdot\Omega_{8,k/8}^2z_{3}
+a\cdot\Omega_{8,k/8}^3z_{4}
-\Omega_{8,k/8}^4z_{5}
-b\cdot\Omega_{8,k/8}^5z_{6}
-i\cdot\Omega_{8,k/8}^6z_{7}
-a\cdot\Omega_{8,k/8}^7z_{8}\hfill\\
z_{1}
-\Omega_{8,k/8}z_{2}
+\Omega_{8,k/8}^2z_{3}
-\Omega_{8,k/8}^3z_{4}
+\Omega_{8,k/8}^4z_{5}
-\Omega_{8,k/8}^5z_{6}
+\Omega_{8,k/8}^6z_{7}
-\Omega_{8,k/8}^7z_{8}\hfill\\
z_{1}
-a\cdot\Omega_{8,k/8}z_{2}
-i\cdot\Omega_{8,k/8}^2z_{3}
-b\cdot\Omega_{8,k/8}^3z_{4}
-\Omega_{8,k/8}^4z_{5}
+a\cdot\Omega_{8,k/8}^5z_{6}
+i\cdot\Omega_{8,k/8}^6z_{7}
+b\cdot\Omega_{8,k/8}^7z_{8}\hfill\\
z_{1}
+i\cdot\Omega_{8,k/8}z_{2}
-\Omega_{8,k/8}^2z_{3}
-i\cdot\Omega_{8,k/8}^3z_{4}
+\Omega_{8,k/8}^4z_{5}
+i\cdot\Omega_{8,k/8}^5z_{6}
-\Omega_{8,k/8}^6z_{7}
-i\cdot\Omega_{8,k/8}^7z_{8}\hfill\\
z_{1}
-b\cdot\Omega_{8,k/8}z_{2}
+i\cdot\Omega_{8,k/8}^2z_{3}
-a\cdot\Omega_{8,k/8}^3z_{4}
-\Omega_{8,k/8}^4z_{5}
+b\cdot\Omega_{8,k/8}^5z_{6}
-i\cdot\Omega_{8,k/8}^6z_{7}
+a\cdot\Omega_{8,k/8}^7z_{8}\hfill
\end{bmatrix}.
\end{align}
As we introduce
\begin{align}
&\tau_1:=z_1+\Omega_{8,k/8}^4z_5,\hspace{1.02cm}
\tau_2:=z_1-\Omega_{8,k/8}^4z_5\\
&\tau_3:=\Omega_{8,k/8}z_2+\Omega_{8,k/8}^5z_6,\;\;
\tau_4:=\Omega_{8,k/8}z_2-\Omega_{8,k/8}^5z_6\nonumber\\
&\tau_5:=\Omega_{8,k/8}^2z_3+\Omega_{8,k/8}^6z_7,\;\;
\tau_6:=\Omega_{8,k/8}^2z_3-\Omega_{8,k/8}^6z_7\nonumber\\
&\tau_7:=\Omega_{8,k/8}^3z_4+\Omega_{8,k/8}^7z_8,\;\;
\tau_8:=\Omega_{8,k/8}^3z_4-\Omega_{8,k/8}^7z_8\nonumber
\end{align}
we obtain
\begin{align}
B_{8,k}x=&\,
\begin{bmatrix}
\tau_1+\tau_3+\tau_5+\tau_7\hfill\\
\tau_2+a\cdot\tau_4-i\cdot\tau_6+b\cdot\tau_8\hfill\\
\tau_1-i\cdot\tau_3-\tau_5+i\cdot\tau_7\hfill\\
\tau_2+b\cdot\tau_4+i\cdot\tau_6+a\cdot\tau_8\hfill\\
\tau_1-\tau_3+\tau_5-\tau_7\hfill\\
\tau_2-a\cdot\tau_4-i\cdot\tau_6-b\cdot\tau_8\hfill\\
\tau_1+i\cdot\tau_3-\tau_5-i\cdot\tau_7\hfill\\
\tau_2-b\cdot\tau_4+i\cdot\tau_6-a\cdot\tau_8\hfill
\end{bmatrix}.
\end{align}
This leads to Algo.~\ref{alg:butterfly-kernel-radix-8} for the computation of $x\mapsto A_{8,n}x$.
\begin{algorithm}[H]
\caption{Butterfly kernel of radix-8}
\label{alg:butterfly-kernel-radix-8}
\begin{algorithmic}[1]
\REQUIRE{$x\in\mathbb{C}^n$, $n=8^t$}
\ENSURE{$x:=A_{8,n}x$}
\FOR{$q=1,\dots,t$}
\STATE{$k:=8^q$}
\STATE{$n_b:=n/k$}
\STATE{$\ell:=k/8$}
\FOR{$b=1,\dots,n_b$}
\FOR{$j=1\dots,\ell$}
\STATE{$z_1:=x_{(b-1)k+j}$}
\STATE{$z_2:=\omega_k^{j-1}x_{(b-1)k+\ell+j}$}
\STATE{$z_3:=\omega_k^{2(j-1)}x_{(b-1)k+2\ell+j}$}
\STATE{$z_4:=\omega_k^{3(j-1)}x_{(b-1)k+3\ell+j}$}
\STATE{$z_5:=\omega_k^{4(j-1)}x_{(b-1)k+4\ell+j}$}
\STATE{$z_6:=\omega_k^{5(j-1)}x_{(b-1)k+5\ell+j}$}
\STATE{$z_7:=\omega_k^{6(j-1)}x_{(b-1)k+6\ell+j}$}
\STATE{$z_8:=\omega_k^{7(j-1)}x_{(b-1)k+7\ell+j}$}
\STATE{$\tau_1:=z_1+z_5$, $\tau_2:=z_1-z_5$}
\STATE{$\tau_3:=z_2+z_6$, $\tau_4:=z_2-z_6$}
\STATE{$\tau_5:=z_3+z_7$, $\tau_6:=z_3-z_7$}
\STATE{$\tau_7:=z_4+z_8$, $\tau_8:=z_4-z_8$}
\STATE{$x_{(b-1)k+j}:=\tau_1+\tau_3+\tau_5+\tau_7$}
\STATE{$x_{(b-1)k+\ell+j}:=\tau_2+a\cdot\tau_4-i\cdot\tau_6+b\cdot\tau_8$}
\STATE{$x_{(b-1)k+2\ell+j}:=\tau_1-i\cdot\tau_3-\tau_5+i\cdot\tau_7$}
\STATE{$x_{(b-1)k+3\ell+j}:=\tau_2+b\cdot\tau_4+i\cdot\tau_6+a\cdot\tau_8$}
\STATE{$x_{(b-1)k+4\ell+j}:=\tau_1-\tau_3+\tau_5-\tau_7$}
\STATE{$x_{(b-1)k+5\ell+j}:=\tau_2-a\cdot\tau_4-i\cdot\tau_6-b\cdot\tau_8$}
\STATE{$x_{(b-1)k+6\ell+j}:=\tau_1+i\cdot\tau_3-\tau_5-i\cdot\tau_7$}
\STATE{$x_{(b-1)k+7\ell+j}:=\tau_2-b\cdot\tau_4+i\cdot\tau_6-a\cdot\tau_8$}
\ENDFOR
\ENDFOR
\ENDFOR
\RETURN{$x$}
\end{algorithmic}
\end{algorithm}
The radix-8 conjugate butterfly kernel relies on the expression
\begin{align}
\hspace{-.4cm}
\overline{B_{8,k}}x=&\,
\begin{bmatrix}
z_{1}
+\overline{\Omega_{8,k/8}}z_{2}
+\overline{\Omega_{8,k/8}^2}z_{3}
+\overline{\Omega_{8,k/8}^3}z_{4}
+\overline{\Omega_{8,k/8}^4}z_{5}
+\overline{\Omega_{8,k/8}^5}z_{6}
+\overline{\Omega_{8,k/8}^6}z_{7}
+\overline{\Omega_{8,k/8}^7}z_{8}\hfill\\
z_{1}
+\overline{a}\cdot\overline{\Omega_{8,k/8}}z_{2}
+i\cdot\overline{\Omega_{8,k/8}^2}z_{3}
+\overline{b}\cdot\overline{\Omega_{8,k/8}^3}z_{4}
-\overline{\Omega_{8,k/8}^4}z_{5}
-\overline{a}\cdot\overline{\Omega_{8,k/8}^5}z_{6}
-i\cdot\overline{\Omega_{8,k/8}^6}z_{7}
-\overline{b}\cdot\overline{\Omega_{8,k/8}^7}z_{8}\hfill\\
z_{1}
+i\cdot\overline{\Omega_{8,k/8}}z_{2}
-\overline{\Omega_{8,k/8}^2}z_{3}
-i\cdot\overline{\Omega_{8,k/8}^3}z_{4}
+\overline{\Omega_{8,k/8}^4}z_{5}
+i\cdot\overline{\Omega_{8,k/8}^5}z_{6}
-\overline{\Omega_{8,k/8}^6}z_{7}
-i\cdot\overline{\Omega_{8,k/8}^7}z_{8}\hfill\\
z_{1}
+\overline{b}\cdot\overline{\Omega_{8,k/8}}z_{2}
-i\cdot\overline{\Omega_{8,k/8}^2}z_{3}
+\overline{a}\cdot\overline{\Omega_{8,k/8}^3}z_{4}
-\overline{\Omega_{8,k/8}^4}z_{5}
-\overline{b}\cdot\overline{\Omega_{8,k/8}^5}z_{6}
+i\cdot\overline{\Omega_{8,k/8}^6}z_{7}
-\overline{a}\cdot\overline{\Omega_{8,k/8}^7}z_{8}\hfill\\
z_{1}
-\overline{\Omega_{8,k/8}}z_{2}
+\overline{\Omega_{8,k/8}^2}z_{3}
-\overline{\Omega_{8,k/8}^3}z_{4}
+\overline{\Omega_{8,k/8}^4}z_{5}
-\overline{\Omega_{8,k/8}^5}z_{6}
+\overline{\Omega_{8,k/8}^6}z_{7}
-\overline{\Omega_{8,k/8}^7}z_{8}\hfill\\
z_{1}
-\overline{a}\cdot\overline{\Omega_{8,k/8}}z_{2}
+i\cdot\overline{\Omega_{8,k/8}^2}z_{3}
-\overline{b}\cdot\overline{\Omega_{8,k/8}^3}z_{4}
-\overline{\Omega_{8,k/8}^4}z_{5}
+\overline{a}\cdot\overline{\Omega_{8,k/8}^5}z_{6}
-i\cdot\overline{\Omega_{8,k/8}^6}z_{7}
+\overline{b}\cdot\overline{\Omega_{8,k/8}^7}z_{8}\hfill\\
z_{1}
-i\cdot\overline{\Omega_{8,k/8}}z_{2}
-\overline{\Omega_{8,k/8}^2}z_{3}
+i\cdot\overline{\Omega_{8,k/8}^3}z_{4}
+\overline{\Omega_{8,k/8}^4}z_{5}
-i\cdot\overline{\Omega_{8,k/8}^5}z_{6}
-\overline{\Omega_{8,k/8}^6}z_{7}
+i\cdot\overline{\Omega_{8,k/8}^7}z_{8}\hfill\\
z_{1}
-\overline{b}\cdot\overline{\Omega_{8,k/8}}z_{2}
-i\cdot\overline{\Omega_{8,k/8}^2}z_{3}
-\overline{a}\cdot\overline{\Omega_{8,k/8}^3}z_{4}
-\overline{\Omega_{8,k/8}^4}z_{5}
+\overline{b}\cdot\overline{\Omega_{8,k/8}^5}z_{6}
+i\cdot\overline{\Omega_{8,k/8}^6}z_{7}
+\overline{a}\cdot\overline{\Omega_{8,k/8}^7}z_{8}\hfill
\end{bmatrix}
\end{align}
so that, as we introduce
\begin{align}
&\tau_1:=z_1+\overline{\Omega_{8,k/8}^4}z_5,\hspace{1.02cm}
\tau_2:=z_1-\overline{\Omega_{8,k/8}^4}z_5\\
&\tau_3:=\overline{\Omega_{8,k/8}}z_2+\overline{\Omega_{8,k/8}^5}z_6,\;\;
\tau_4:=\overline{\Omega_{8,k/8}}z_2-\overline{\Omega_{8,k/8}^5}z_6\nonumber\\
&\tau_5:=\overline{\Omega_{8,k/8}^2}z_3+\overline{\Omega_{8,k/8}^6}z_7,\;\;
\tau_6:=\overline{\Omega_{8,k/8}^2}z_3-\overline{\Omega_{8,k/8}^6}z_7\nonumber\\
&\tau_7:=\overline{\Omega_{8,k/8}^3}z_4+\overline{\Omega_{8,k/8}^7}z_8,\;\;
\tau_8:=\overline{\Omega_{8,k/8}^3}z_4-\overline{\Omega_{8,k/8}^7}z_8\nonumber
\end{align}
we obtain
\begin{align}
\overline{B_{8,k}}x=&\,
\begin{bmatrix}
\tau_1+\tau_3+\tau_5+\tau_7\hfill\\
\tau_2+\overline{a}\cdot\tau_4+i\cdot\tau_6+\overline{b}\cdot\tau_8\hfill\\
\tau_1+i\cdot\tau_3-\tau_5-i\cdot\tau_7\hfill\\
\tau_2+\overline{b}\cdot\tau_4-i\cdot\tau_6+\overline{a}\cdot\tau_8\hfill\\
\tau_1-\tau_3+\tau_5-\tau_7\hfill\\
\tau_2-\overline{a}\cdot\tau_4+i\cdot\tau_6-\overline{b}\cdot\tau_8\hfill\\
\tau_1-i\cdot\tau_3-\tau_5+i\cdot\tau_7\hfill\\
\tau_2-\overline{b}\cdot\tau_4-i\cdot\tau_6-\overline{a}\cdot\tau_8\hfill
\end{bmatrix}.
\end{align}
This leads to Algo.~\ref{alg:conjugate-butterfly-kernel-radix-8} for the computation of $x\mapsto \overline{A_{8,n}}x$.
\begin{algorithm}[H]
\caption{Conjugate butterfly kernel of radix-8}
\label{alg:conjugate-butterfly-kernel-radix-8}
\begin{algorithmic}[1]
\REQUIRE{$x\in\mathbb{C}^n$, $n=8^t$}
\ENSURE{$x:=\overline{A_{8,n}}x$}
\FOR{$q=1,\dots,t$}
\STATE{$k:=8^q$}
\STATE{$n_b:=n/k$}
\STATE{$\ell:=k/8$}
\FOR{$b=1,\dots,n_b$}
\FOR{$j=1\dots,\ell$}
\STATE{$z_1:=x_{(b-1)k+j}$}
\STATE{$z_2:=\overline{\omega_k^{j-1}}x_{(b-1)k+\ell+j}$}
\STATE{$z_3:=\overline{\omega_k^{2(j-1)}}x_{(b-1)k+2\ell+j}$}
\STATE{$z_4:=\overline{\omega_k^{3(j-1)}}x_{(b-1)k+3\ell+j}$}
\STATE{$z_5:=\overline{\omega_k^{4(j-1)}}x_{(b-1)k+4\ell+j}$}
\STATE{$z_6:=\overline{\omega_k^{5(j-1)}}x_{(b-1)k+5\ell+j}$}
\STATE{$z_7:=\overline{\omega_k^{6(j-1)}}x_{(b-1)k+6\ell+j}$}
\STATE{$z_8:=\overline{\omega_k^{7(j-1)}}x_{(b-1)k+7\ell+j}$}
\STATE{$\tau_1:=z_1+z_5$, $\tau_2:=z_1-z_5$}
\STATE{$\tau_3:=z_2+z_6$, $\tau_4:=z_2-z_6$}
\STATE{$\tau_5:=z_3+z_7$, $\tau_6:=z_3-z_7$}
\STATE{$\tau_7:=z_4+z_8$, $\tau_8:=z_4-z_8$}
\STATE{$x_{(b-1)k+j}:=\tau_1+\tau_3+\tau_5+\tau_7$}
\STATE{$x_{(b-1)k+\ell+j}:=\tau_2+\overline{a}\cdot\tau_4+i\cdot\tau_6+\overline{b}\cdot\tau_8$}
\STATE{$x_{(b-1)k+2\ell+j}:=\tau_1+i\cdot\tau_3-\tau_5-i\cdot\tau_7$}
\STATE{$x_{(b-1)k+3\ell+j}:=\tau_2+\overline{b}\cdot\tau_4-i\cdot\tau_6+\overline{a}\cdot\tau_8$}
\STATE{$x_{(b-1)k+4\ell+j}:=\tau_1-\tau_3+\tau_5-\tau_7$}
\STATE{$x_{(b-1)k+5\ell+j}:=\tau_2-\overline{a}\cdot\tau_4+i\cdot\tau_6-\overline{b}\cdot\tau_8$}
\STATE{$x_{(b-1)k+6\ell+j}:=\tau_1-i\cdot\tau_3-\tau_5+i\cdot\tau_7$}
\STATE{$x_{(b-1)k+7\ell+j}:=\tau_2-\overline{b}\cdot\tau_4-i\cdot\tau_6-\overline{a}\cdot\tau_8$}
\ENDFOR
\ENDFOR
\ENDFOR
\RETURN{$x$}
\end{algorithmic}
\end{algorithm}
On the other hand, we have
\begin{align}
\hspace{-.8cm}
B_{8,k}^Tx=&\,
\begin{bmatrix}
z_{1}+z_{2}+z_{3}+z_{4}+z_{5}+z_{6}+z_{7}+z_{8}\hfill\\
\Omega_{8,k/8}z_{1}
+a\cdot\Omega_{8,k/8}z_{2}
-i\cdot\Omega_{8,k/8}z_{3}
+b\cdot\Omega_{8,k/8}z_{4}
-\Omega_{8,k/8}z_{5}
-a\cdot\Omega_{8,k/8}z_{6}
+i\cdot\Omega_{8,k/8}z_{7}
-b\cdot\Omega_{8,k/8}z_{8}\hfill\\
\Omega_{8,k/8}^2z_{1}
-i\cdot\Omega_{8,k/8}^2z_{2}
-\Omega_{8,k/8}^2z_{3}
+i\cdot\Omega_{8,k/8}^2z_{4}
+\Omega_{8,k/8}^2z_{5}
-i\cdot\Omega_{8,k/8}^2z_{6}
-\Omega_{8,k/8}^2z_{7}
+i\cdot\Omega_{8,k/8}^2z_{8}\hfill\\
\Omega_{8,k/8}^3z_{1}
+b\cdot\Omega_{8,k/8}^3z_{2}
+i\cdot\Omega_{8,k/8}^3z_{3}
+a\cdot\Omega_{8,k/8}^3z_{4}
-\Omega_{8,k/8}^3z_{5}
-b\cdot\Omega_{8,k/8}^3z_{6}
-i\cdot\Omega_{8,k/8}^3z_{7}
-a\cdot\Omega_{8,k/8}^3z_{8}\hfill\\
\Omega_{8,k/8}^4z_{1}
-\Omega_{8,k/8}^4z_{2}
+\Omega_{8,k/8}^4z_{3}
-\Omega_{8,k/8}^4z_{4}
+\Omega_{8,k/8}^4z_{5}
-\Omega_{8,k/8}^4z_{6}
+\Omega_{8,k/8}^4z_{7}
-\Omega_{8,k/8}^4z_{8}\hfill\\
\Omega_{8,k/8}^5z_{1}
-a\cdot\Omega_{8,k/8}^5z_{2}
-i\cdot\Omega_{8,k/8}^5z_{3}
-b\cdot\Omega_{8,k/8}^5z_{4}
-\Omega_{8,k/8}^5z_{5}
+a\cdot\Omega_{8,k/8}^5z_{6}
+i\cdot\Omega_{8,k/8}^5z_{7}
+b\cdot\Omega_{8,k/8}^5z_{8}\hfill\\
\Omega_{8,k/8}^6z_{1}
+i\cdot\Omega_{8,k/8}^6z_{2}
-\Omega_{8,k/8}^6z_{3}
-i\cdot\Omega_{8,k/8}^6z_{4}
+\Omega_{8,k/8}^6z_{5}
+i\cdot\Omega_{8,k/8}^6z_{6}
-\Omega_{8,k/8}^6z_{7}
-i\cdot\Omega_{8,k/8}^6z_{8}\hfill\\
\Omega_{8,k/8}^7z_{1}
-b\cdot\Omega_{8,k/8}^7z_{2}
+i\cdot\Omega_{8,k/8}^7z_{3}
-a\cdot\Omega_{8,k/8}^7z_{4}
-\Omega_{8,k/8}^7z_{5}
+b\cdot\Omega_{8,k/8}^7z_{6}
-i\cdot\Omega_{8,k/8}^7z_{7}
+a\cdot\Omega_{8,k/8}^7z_{8}\hfill
\end{bmatrix}.
\end{align}
As we introduce
\begin{align}
\tau_1:=z_1+z_5,\;\;
\tau_2:=\Omega_{8,k/8}(z_1-z_5)\\
\tau_3:=z_2+z_6,\;\;
\tau_4:=\Omega_{8,k/8}(z_2-z_6)\nonumber\\
\tau_5:=z_3+z_7,\;\;
\tau_6:=\Omega_{8,k/8}(z_3-z_7)\nonumber\\
\tau_7:=z_4+z_8,\;\;
\tau_8:=\Omega_{8,k/8}(z_4-z_8)\nonumber
\end{align}
we obtain
\begin{align}
B_{8,k}^Tx=&\,
\begin{bmatrix}
\tau_1+\tau_3+\tau_5+\tau_7\hfill\\
\tau_2+a\cdot\tau_4-i\cdot\tau_6+b\cdot\tau_8\hfill\\
\Omega_{8,k/8}^2(\tau_1-i\cdot\tau_3-\tau_5+i\cdot\tau_7)\hfill\\
\Omega_{8,k/8}^2(\tau_2+b\cdot\tau_4+i\cdot\tau_6+a\cdot\tau_8)\hfill\\
\Omega_{8,k/8}^4(\tau_1-\tau_3+\tau_5-\tau_7)\hfill\\
\Omega_{8,k/8}^4(\tau_2-a\cdot\tau_4-i\cdot\tau_6-b\cdot\tau_8)\hfill\\
\Omega_{8,k/8}^6(\tau_1+i\cdot\tau_3-\tau_5-i\cdot\tau_7)\hfill\\
\Omega_{8,k/8}^6(\tau_2-b\cdot\tau_4+i\cdot\tau_6-a\cdot\tau_8)\hfill
\end{bmatrix}.
\end{align}
This leads to Algo.~\ref{alg:transposed-butterfly-kernel-radix-8} for the computation of $x\mapsto A_{8,n}^Tx$.
\begin{algorithm}[H]
\caption{Transposed butterfly kernel of radix-8}
\label{alg:transposed-butterfly-kernel-radix-8}
\begin{algorithmic}[1]
\REQUIRE{$x\in\mathbb{C}^n$, $n=8^t$}
\ENSURE{$x:=A_{8,n}^Tx$}
\FOR{$q=t,\dots,1$}
\STATE{$k:=8^q$}
\STATE{$n_b:=n/k$}
\STATE{$\ell:=k/8$}
\FOR{$b=1,\dots,n_b$}
\FOR{$j=1\dots,\ell$}
\STATE{$z_1:=x_{(b-1)k+j}$}
\STATE{$z_2:=x_{(b-1)k+\ell+j}$}
\STATE{$z_3:=x_{(b-1)k+2\ell+j}$}
\STATE{$z_4:=x_{(b-1)k+3\ell+j}$}
\STATE{$z_5:=x_{(b-1)k+4\ell+j}$}
\STATE{$z_6:=x_{(b-1)k+5\ell+j}$}
\STATE{$z_7:=x_{(b-1)k+6\ell+j}$}
\STATE{$z_8:=x_{(b-1)k+7\ell+j}$}
\STATE{$\tau_1:=z_1+z_5$, $\tau_2:=\omega_k^{j-1}(z_1-z_5)$}
\STATE{$\tau_3:=z_2+z_6$, $\tau_4:=\omega_k^{j-1}(z_2-z_6)$}
\STATE{$\tau_5:=z_3+z_7$, $\tau_6:=\omega_k^{j-1}(z_3-z_7)$}
\STATE{$\tau_7:=z_4+z_8$, $\tau_8:=\omega_k^{j-1}(z_4-z_8)$}
\STATE{$x_{(b-1)k+j}:=\tau_1+\tau_3+\tau_5+\tau_7$}
\STATE{$x_{(b-1)k+\ell+j}:=\tau_2+a\cdot\tau_4-i\cdot\tau_6+b\cdot\tau_8$}
\STATE{$x_{(b-1)k+2\ell+j}:=\omega_{k}^{2(j-1)}(\tau_1-i\cdot\tau_3-\tau_5+i\cdot\tau_7)$}
\STATE{$x_{(b-1)k+3\ell+j}:=\omega_{k}^{2(j-1)}(\tau_2+b\cdot\tau_4+i\cdot\tau_6+a\cdot\tau_8)$}
\STATE{$x_{(b-1)k+4\ell+j}:=\omega_{k}^{4(j-1)}(\tau_1-\tau_3+\tau_5-\tau_7)$}
\STATE{$x_{(b-1)k+5\ell+j}:=\omega_{k}^{4(j-1)}(\tau_2-a\cdot\tau_4-i\cdot\tau_6-b\cdot\tau_8)$}
\STATE{$x_{(b-1)k+6\ell+j}:=\omega_{k}^{6(j-1)}(\tau_1+i\cdot\tau_3-\tau_5-i\cdot\tau_7)$}
\STATE{$x_{(b-1)k+7\ell+j}:=\omega_{k}^{6(j-1)}(\tau_2-b\cdot\tau_4+i\cdot\tau_6-a\cdot\tau_8)$}
\ENDFOR
\ENDFOR
\ENDFOR
\RETURN{$x$}
\end{algorithmic}
\end{algorithm}

\section{2D isotropic numerical experiments}\label{app-2d}

\begin{table}
\caption{Timing of 2D index-reversal permutation, unordered and full forward FFT --- Computing system \#1}
\label{tab:2d-fft}
\begin{tabular}{cccccc}
\toprule
Procedure     & 
$n_1=n_2=2^{10}$& 
$n_1=n_2=2^{11}$& 
$n_1=n_2=2^{12}$& 
$n_1=n_2=2^{13}$& 
$n_1=n_2=2^{14}$\\
\midrule
\multicolumn{6}{c}{radix-2}\\
\shortstack{Index-reversal\\permutation} 
              &  0.002~s & 0.007~s & 0.092~s & 0.423~s &  1.738~s\\
Unordered FFT &  0.027~s & 0.192~s & 1.005~s & 4.409~s & 19.335~s\\
FFT           &  0.029~s & 0.197~s & 1.100~s & 4.827~s & 21.066~s\\
\midrule
\multicolumn{6}{c}{radix-4}\\
\shortstack{Index-reversal\\permutation}   
              &  0.001~s & -       & 0.096~s & -       &  1.675~s\\
Unordered FFT &  0.027~s & -       & 1.006~s & -       & 19.248~s\\
FFT           &  0.029~s & -       & 1.101~s & -       & 20.903~s\\
\midrule
FFT (FFTW)    &  0.025~s & 0.142~s & 0.993~s & 5.660~s & 32.805~s\\
\bottomrule
\end{tabular}
\end{table}

\begin{table}
\caption{Timing of 2D index-reversal permutation, unordered and full backward FFT --- Computing system \#1}
\label{tab:2d-ifft}
\begin{tabular}{cccccc}
\toprule
Procedure     & 
$n_1=n_2=2^{10}$& 
$n_1=n_2=2^{11}$& 
$n_1=n_2=2^{12}$& 
$n_1=n_2=2^{13}$& 
$n_1=n_2=2^{14}$\\
\midrule
\multicolumn{6}{c}{radix-2}\\
\shortstack{Index-reversal\\permutation}
              &  0.002~s & 0.007~s & 0.091~s & 0.427~s &  1.740~s\\
Unordered FFT &  0.027~s & 0.196~s & 1.029~s & 4.447~s & 19.339~s\\
FFT           &  0.029~s & 0.203~s & 1.122~s & 4.871~s & 21.088~s\\
\midrule
\multicolumn{6}{c}{radix-4}\\
\shortstack{Index-reversal\\permutation}
              &  0.001~s & -       & 0.096~s & -       &  1.675~s\\
Unordered FFT &  0.028~s & -       & 1.032~s & -       & 19.415~s\\
FFT           &  0.029~s & -       & 1.145~s & -       & 21.018~s\\
\midrule
FFT (FFTW)    &  0.021~s & 0.104~s & 0.083~s & 5.268~s & 32.503~s\\
\bottomrule
\end{tabular}
\end{table}

\begin{table}
\caption{Timing of 2D index-reversal permutation, unordered and full forward FFT --- Computing system \#2}
\label{tab:2d-fft-2}
\begin{tabular}{cccccc}
\toprule
Procedure     & 
$n_1=n_2=2^{10}$& 
$n_1=n_2=2^{11}$& 
$n_1=n_2=2^{12}$& 
$n_1=n_2=2^{13}$& 
$n_1=n_2=2^{14}$\\
\midrule
\multicolumn{6}{c}{radix-2}\\
\shortstack{Index-reversal\\permutation} 
              &  0.002~s & 0.008~s & 0.071~s & 0.346~s &  1.465~s\\
Unordered FFT &  0.035~s & 0.141~s & 1.138~s & 5.762~s & 23.968~s\\
FFT           &  0.037~s & 0.149~s & 1.213~s & 6.109~s & 25.465~s\\
\midrule
\multicolumn{6}{c}{radix-4}\\
\shortstack{Index-reversal\\permutation}   
              &  0.002~s & -       & 0.066~s & -       &  1.408~s\\
Unordered FFT &  0.034~s & -       & 1.100~s & -       & 23.565~s\\
FFT           &  0.035~s & -       & 1.185~s & -       & 24.942~s\\
\midrule
FFT (FFTW)    &  0.031~s & 0.162~s & 1.647~s & 7.258~s & 41.890~s\\
\bottomrule
\end{tabular}
\end{table}

\begin{table}
\caption{Timing of 2D index-reversal permutation, unordered and full backward FFT         --- Computing system \#2}
\label{tab:2d-ifft-2}
\begin{tabular}{cccccc}
\toprule
Procedure     & 
$n_1=n_2=2^{10}$& 
$n_1=n_2=2^{11}$& 
$n_1=n_2=2^{12}$& 
$n_1=n_2=2^{13}$& 
$n_1=n_2=2^{14}$\\
\midrule
\multicolumn{6}{c}{radix-2}\\
\shortstack{Index-reversal\\permutation}
              &  0.002~s & 0.008~s & 0.072~s & 0.346~s &  1.496~s\\
Unordered FFT &  0.035~s & 0.144~s & 1.153~s & 5.826~s & 24.031~s\\
FFT           &  0.037~s & 0.152~s & 1.224~s & 6.172~s & 25.526~s\\
\midrule
\multicolumn{6}{c}{radix-4}\\
\shortstack{Index-reversal\\permutation}   
              &  0.002~s & -       & 0.066~s & -       &  1.411~s\\
Unordered FFT &  0.034~s & -       & 1.137~s & -       & 23.762~s\\
FFT           &  0.036~s & -       & 1.205~s & -       & 25.171~s\\
\midrule
FFT (FFTW)    &  0.032~s & 0.164~s & 1.498~s & 7.198~s & 42.803~s\\
\bottomrule
\end{tabular}
\end{table}

\section{3D isotropic numerical experiments}\label{app-3d}

\begin{table}
\caption{Timing of 3D index-reversal permutation, unordered and full forward FFT --- Computing system \#1}
\label{tab:3d-fft}
\begin{tabular}{cccc}
\toprule
Procedure     & 
$n_1=n_2=n_3=2^{7}$& 
$n_1=n_2=n_3=2^{8}$& 
$n_1=n_2=n_3=2^{9}$\\
\midrule
\multicolumn{4}{c}{radix-2}\\
\shortstack{Index-reversal\\permutation} 
              &  0.002~s & 0.022~s & 0.194~s\\
Unordered FFT &  0.047~s & 0.544~s & 5.841~s\\
FFT           &  0.050~s & 0.561~s & 6.048~s\\
\midrule
\multicolumn{4}{c}{radix-4}\\
\shortstack{Index-reversal\\permutation}   
              &  -       & 0.021~s & -      \\
Unordered FFT &  -       & 0.542~s & -      \\
FFT           &  -       & 0.560~s & -      \\
\midrule
FFT (FFTW)    &  0.013~s & 0.498~s & 5.702~s\\
\bottomrule
\end{tabular}
\end{table}

\begin{table}
\caption{Timing of 3D index-reversal permutation, unordered and full backward FFT --- Computing system \#1}
\label{tab:3d-ifft}
\begin{tabular}{cccc}
\toprule
Procedure     & 
$n_1=n_2=n_3=2^{7}$& 
$n_1=n_2=n_3=2^{8}$& 
$n_1=n_2=n_3=2^{9}$\\
\midrule
\multicolumn{4}{c}{radix-2}\\
\shortstack{Index-reversal\\permutation}
              &  0.002~s & 0.022~s & 0.197~s\\
Unordered FFT &  0.048~s & 0.558~s & 5.935~s\\
FFT           &  0.051~s & 0.581~s & 6.138~s\\
\midrule
\multicolumn{4}{c}{radix-4}\\
\shortstack{Index-reversal\\permutation}
              &  -       & 0.020~s & -      \\
Unordered FFT &  -       & 0.550~s & -      \\
FFT           &  -       & 0.571~s & -      \\
\midrule
FFT (FFTW)    &  0.014~s & 0.491~s & 5.712~s\\
\bottomrule
\end{tabular}
\end{table}

\begin{table}
\caption{Timing of 3D index-reversal permutation, unordered and full forward FFT --- Computing system \#2}
\label{tab:3d-fft-2}
\begin{tabular}{cccc}
\toprule
Procedure     & 
$n_1=n_2=n_3=2^{7}$& 
$n_1=n_2=n_3=2^{8}$& 
$n_1=n_2=n_3=2^{9}$\\
\midrule
\multicolumn{4}{c}{radix-2}\\
\shortstack{Index-reversal\\permutation} 
              &  0.002~s & 0.029~s &  0.279~s\\
Unordered FFT &  0.081~s & 0.833~s & 10.458~s\\
FFT           &  0.083~s & 0.862~s & 10.733~s\\
\midrule
\multicolumn{4}{c}{radix-4}\\
\shortstack{Index-reversal\\permutation}   
              &  -       & 0.030~s & -      \\
Unordered FFT &  -       & 0.788~s & -      \\
FFT           &  -       & 0.820~s & -      \\
\midrule
FFT (FFTW)    &  0.015~s & 1.003~s & 11.232~s\\
\bottomrule
\end{tabular}
\end{table}

\begin{table}
\caption{Timing of 3D index-reversal permutation, unordered and full backward FFT --- Computing system \#2}
\label{tab:3d-ifft-2}
\begin{tabular}{cccc}
\toprule
Procedure     & 
$n_1=n_2=n_3=2^{7}$& 
$n_1=n_2=n_3=2^{8}$& 
$n_1=n_2=n_3=2^{9}$\\
\midrule
\multicolumn{4}{c}{radix-2}\\
\shortstack{Index-reversal\\permutation}
              &  0.002~s & 0.030~s &  0.284~s\\
Unordered FFT &  0.082~s & 0.847~s & 10.631~s\\
FFT           &  0.084~s & 0.879~s & 10.909~s\\
\midrule
\multicolumn{4}{c}{radix-4}\\
\shortstack{Index-reversal\\permutation}
              &  -       & 0.031~s & -      \\
Unordered FFT &  -       & 0.804~s & -      \\
FFT           &  -       & 0.835~s & -      \\
\midrule
FFT (FFTW)    &  0.016~s & 1.021~s & 11.524~s\\
\bottomrule
\end{tabular}
\end{table}

\end{document}